\documentclass[11pt]{article}
\usepackage{cite}
\usepackage{amsmath,amssymb,latexsym}
\usepackage{multirow}

\begin{document}
\newcommand{\qed}{\hphantom{.}\hfill $\Box$\medbreak}
\newcommand{\proof}{\noindent{\bf Proof \ }}
\renewcommand{\theequation}{\thesection.\arabic{equation}}
\newtheorem{Theorem}{Theorem}[section]
\newtheorem{Lemma}[Theorem]{Lemma}
\newtheorem{Corollary}[Theorem]{Corollary}
\newtheorem{Proposition}[Theorem]{Proposition}
\newtheorem{Remark}[Theorem]{Remark}
\newtheorem{Example}[Theorem]{Example}
\newtheorem{Definition}[Theorem]{Definition}
\newtheorem{Construction}[Theorem]{Construction}

\begin{center}
{\Large\bf Optimal optical orthogonal signature pattern codes with weight three and cross-correlation constraint one
%\footnote{Supported by NSFC under Grant 11601472, and the Yunnan Applied Basic Research Project of China under Grant 2016FD005 (R. Pan), NSFC under Grant $11871095$, and Fundamental Research Funds for the Central Universities under Grant $2016$JBZ$012$ (T. Feng), NSFC under Grant $11401582$, NSFHB under Grant A2015507019, and Cultivation Project of National Natural Science Foundation of China People's Police University under Grant ZKJJPY201703 (L. Wang), NSFC under Grant 11771227, and Zhejiang Provincial Natural Science Foundation of China under Grant LY17A010008 (X. Wang).}
\footnote{Supported by NSFC under Grant 11601472, and the Yunnan Applied Basic Research Project of China under Grant 2016FD005 (R. Pan), NSFC under Grant $11871095$ (T. Feng), NSFC under Grant 11771227 (X. Wang).}
}

\vskip12pt

Rong Pan$^{a}$, Tao Feng$^{b}$, Lidong Wang$^{c}$, Xiaomiao Wang$^{d}$\\[2ex]
{\footnotesize $^a$Department of Mathematics, Yunnan University, Kunming 650504, P. R. China}\\
{\footnotesize $^b$Department of Mathematics, Beijing Jiaotong University, Beijing 100044, P. R. China}\\
{\footnotesize $^c$Department of Basic Courses, China People's Police University, Langfang 065000, P. R. China}\\
{\footnotesize $^d$Department of Mathematics, Ningbo University, Ningbo 315211, P. R. China}\\
{\footnotesize p$_{-}$rita@163.com, tfeng@bjtu.edu.cn, lidongwang@aliyun.com, wangxiaomiao@nbu.edu.cn}
\vskip12pt

\end{center}

\noindent {\bf Abstract:}
Optical orthogonal signature pattern codes (OOSPCs) have attracted wide attention as signature patterns of spatial optical code division multiple access networks. In this paper, an improved upper bound on the size of an $(m,n,3,\lambda_a,1)$-OOSPC with $\lambda_a=2,3$ is established. The exact number of codewords of an optimal $(m,n,3,\lambda_a,1)$-OOSPC is determined for any positive integers $m,n\equiv2\ ({\rm mod }\ 4)$ and $\lambda_a\in\{2,3\}$.



%%%%%%%%%%%%%%%%%%%%%%%%%%%%%%%%%%%%%%%%%%%%%%%%%%%%%%%%%%%%%%%%%%%%%%%%%

\section{Introduction}

An optical orthogonal signature pattern code is a family of $(0,1)$-matrices with good auto- and cross-correlation. Its study has been motivated by an application in an optical code division
multiple access (OCDMA) network for image transmission, called a spatial OCDMA network.
%The spatial OCDMA network has promoted the development of high-speed multiple access network applications, especially image applications such as supercomputer visualizations, medical image access, and distribution and digital video broadcasting.
Compared with the traditional OCDMA, the spatial OCDMA provides higher throughput (cf. \cite{k,KK,KY,GW}).

Denote by $\mathbb{Z}_v$ the additive group of integers modulo $v$. Let $m$, $n$, $k$, $\lambda_a$ and $\lambda_c$ be positive integers. An $(m, n, k, \lambda_a, \lambda_c)$ {\em optical orthogonal signature
pattern code} $($briefly, $(m, n, k, \lambda_a, \lambda_c)$-{\rm OOSPC}$)$
is a family ${\cal C}$ of $m \times n$ $(0,1)$-matrices of Hamming
weight $k$ satisfying the following properties:
\begin{enumerate}
\item[{\rm (1)}] the auto-correlation property: $\sum\limits_{i=0}^{m-1}\sum\limits_{j=0}^{n-1}x_{i,j}x_{i\oplus s,j\widehat{\oplus} t}\leq\lambda_{a}$ for any $(x_{ij})\in {\cal C}$ and any $(s,t)\in \mathbb{Z}_m\times \mathbb{Z}_n\setminus \{(0,0)\};$
\item[{\rm (2)}] the cross-correlation property: $\sum\limits_{i=0}^{m-1}\sum\limits_{j=0}^{n-1}x_{i,j}y_{i\oplus s,j\widehat{\oplus} t}\leq\lambda_{c}$ for any distinct $(x_{ij}), (y_{ij})\in {\cal C}$ and any $(s,t)\in \mathbb{Z}_m\times \mathbb{Z}_n,$
\end{enumerate}
where the additions $\oplus$ and $\widehat{\oplus}$ are, respectively, reduced modulo
$m$ and $n$. When $\lambda_{a}=\lambda_{c}=\lambda$, the notation
$(m, n, k, \lambda_a, \lambda_c)$-OOSPC is briefly written as
$(m, n, k, \lambda)$-OOSPC.

The number of codewords in an OOSPC is called the \emph{size} of the OOSPC. For given positive integers $m$, $n$, $k$, $\lambda_a$ and $\lambda_c$, denote by $\Theta(m,n,k,\lambda_a,\lambda_c)$ the largest possible size among all $(m,n,k,\lambda_a,\lambda_c)$-OOSPCs. An $(m,n,k,\lambda_a,\lambda_c)$-OOSPC with size
$\Theta(m,n,k,\lambda_a,\lambda_c)$ is said to be {\em optimal} (or {\em maximum}).

When $\lambda_{a}=\lambda_{c}=\lambda$, $\Theta(m,n,k,\lambda_a,\lambda_c)$ is simply written as
$\Theta(m,n,k,\lambda)$. Based on the Johnson bound \cite{J} for constant weight codes, an upper bound on $\Theta(m,n,k,\lambda)$ was given below
\begin{eqnarray}\nonumber
\label{eqq1}
\Theta(m,n,k,\lambda)\leq J(mn,k,\lambda),
\end{eqnarray}
\noindent where
$$J(mn,k,\lambda)=\Big\lfloor\frac{1}{k}\Big\lfloor\frac{mn-1}{k-1}\Big\lfloor\frac{mn-2}{k-2}\Big\lfloor\cdots
\Big\lfloor\frac{mn-\lambda}{k-\lambda}\Big\rfloor\cdots\Big\rfloor\Big\rfloor\Big\rfloor\Big\rfloor.
$$
%and $\lfloor x\rfloor$ denotes the largest integer not exceeding $x$.
When $\lambda_a>\lambda_c$, $\Theta(m,n,k,\lambda_a,\lambda_c)$ is upper bounded in \cite{GW} by
\begin{eqnarray}\label{eqq2}
\left\lfloor\frac{\lambda_a(mn-1)(mn-2)\cdots(mn-\lambda_c)}{k(k-1)(k-2)\cdots(k-\lambda_c)}\right\rfloor.
\end{eqnarray}

When $m$ and $n$ are coprime, it has been shown in \cite{GW} that an $(m,n,k,\lambda_a,\lambda_c)$-OOSPC is equivalent to a $1$-dimensional $(mn,k,\lambda_a,\lambda_c)$-optical orthogonal code (OOC). See \cite{JAM,BM,BP,BPW,CSW,FCJ2,FM,GY,Y} and the references therein for more details on OOCs.

When $m$ and $n$ are not coprime, various OOSPCs have been constructed via algebraic and combinatorial methods for the case of $\lambda_a=\lambda_c$ (see \cite{CJL,CJL2,JDWG,PC,PC3,MS,sawa,GW}). We only quote the following result for later use.

\begin{Theorem}\label{thm-pan} {\rm \cite{PC3}}
\begin{eqnarray}\nonumber
\Theta(m,n,3,1)=\left\{
	\begin{array}{ll}
    J(mn,3,1)-1, & \hbox{ \rm{if} $mn\equiv 14,20\ ({\rm mod }\ 24)$,}\\[0.4em]
	&\ \rm{or}\ \hbox{ \rm{if} $mn\equiv 8,16\ ({\rm mod }\ 24)$}\ \rm{and}\ \hbox{$\gcd(m,n,4)=2$},\\[0.4em]
    &\ \rm{or}\ \hbox{ \rm{if} $mn\equiv 2\ ({\rm mod }\ 6)$}\ \rm{and}\ \hbox{$\gcd(m,n,4)=4$};\\[0.4em]
	J(mn,3,1), &\ \rm{otherwise}.\\[0.4em]
	\end{array}
	\right.
	\end{eqnarray}
\end{Theorem}

On the other hand, for the case of $\lambda_a\neq\lambda_c$, very little has been done
on $(m,n,k,\lambda_a,\lambda_c)$-OOSPCs with maximum size.
%To our knowledge, Yang and Kwong \cite{GW} presented three algebraic constructions for determining the value $\Theta(m,m,k,\lambda_a,1)=\lfloor\frac{\lambda_a(m^{2}-1)}{k(k-1)}\rfloor$ for the following parameters:
%\begin{enumerate}
%\item[{\rm (1)}]$(m,k,\lambda_a)=(p,\frac{p-1}{2},\frac{p-3}{4})$
%where $p\equiv 3\ ({\rm mod }\ 4)$ be a prime;
%\item[{\rm (2)}] $(m,k,\lambda_a)=(\frac{q^{t+1}-1}{q-1},\frac{q^{t}-1}{q-1},
%\frac{q^{t-1}-1}{q-1})$ where $q$ be a prime power and $t$ be a positive integer;
%\item[{\rm (3)}] $(m,k,\lambda_a)=(p,k,1)$ where $p\equiv 1\ ({\rm mod }\ k(k-1))$
%be a prime.
%\end{enumerate}
Compared with \eqref{eqq2}, an improved upper bound on $\Theta(m,n,3,2,1)$ was given by Sawa and Kageyama \cite{sawa}. That is
\begin{eqnarray}\label{eqq3}
\Theta(m,n,3,2,1)\leq\left\{
    \begin{array}{ll}
        \frac{mn}{4}, &\hbox{if $mn\equiv0\ ({\rm mod }\ 4)$,}\\[0.4em]
        \left\lfloor\frac{mn-1}{4}\right\rfloor,&\hbox{otherwise.}
    \end{array}
    \right.
\end{eqnarray}
And they proved the following theorem.

\begin{Theorem}\label{thm-sawa} {\rm \cite{sawa}}
\begin{eqnarray}\label{eqqgg}\nonumber
\Theta(m,n,3,2,1)=\left\{
    \begin{array}{ll}
        \frac{mn-1}{4}, &\hbox{\rm{if} $m=n\equiv1\ ({\rm mod }\ 4)$ {\rm is a prime and} $2$ is a primitive root in $\mathbb{Z}_m$};\\
%&\hbox{and $2$ is a primitive root in $\mathbb{Z}_m$,}\\[0.2em]
        \frac{mn-2}{4},&\hbox{\rm{if} $mn\equiv2\ ({\rm mod }\ 4)$.}
    \end{array}
    \right.
\end{eqnarray}
\end{Theorem}

%Recently, there are many studies on multiple-weight OOSPCs, which are generalizations of constant-weight OOSPCs. It was proposed by Kwong and Yang \cite{KY} for OCDMA networks to support multiple quality of services (QoS). The reader may refer to \cite{KY,KY2,ZQ} and the references therein for more detail results.
In Section \ref{2}, we shall give an equivalent combinatorial description of $(m,n,k,\lambda_a,\lambda_c)$-OOSPCs by using set-theoretic notation. Section \ref{331} is devoted to improving Sawa and Kageyama's bound \eqref{eqq3}, especially for the case of $mn\equiv0\ ({\rm mod }\ 4)$. Throughout this paper, let $\xi$ denote the number of subgroups of order $3$ in $\mathbb{Z}_m\times \mathbb{Z}_n$, i.e.,
$$\xi=\left\{
	\begin{array}{ll}
	0,\ \hbox{  \rm{if} $3\nmid mn$\rm{;}}\\
	1,\ \hbox{ \rm{if}  $3\mid mn$ and $\gcd(m,n,3)=1$\rm{;}}\\
	4,\ \hbox{ \rm{if} $\gcd(m,n,3)=3$.}\\
	\end{array}
	\right.
$$
Let
$$\omega=\left\{
	\begin{array}{ll}
	0,\ \hbox{  \rm{if} $\lambda_a=2$;}\\
	\xi,\ \hbox{ \rm{if}  $\lambda_a=3$.}\\
	\end{array}
	\right.
$$
We are to prove the following theorem.
\begin{Theorem}\label{c-main-1}
Let $\lambda_a\in\{2,3\}$. Then $\Theta(m,n,3,\lambda_a,1)\leq$
\begin{eqnarray}\nonumber
\left\{
	\begin{array}{ll}
    \left\lfloor\frac{ mn+2\omega}{4}\right\rfloor,
	& \hbox{$mn\equiv 1,2,3\ ({\rm mod }\ 4)$\rm{;}}\\[0.4em]
	
	\left\lfloor\frac{7mn+16\omega}{32}\right\rfloor,&
	\hbox{$mn\equiv 0\ ({\rm mod }\ 8)$ \rm{and} $\gcd(m,n,2)=1$\rm{;}}\\[0.4em]
	
	\left\lfloor\frac{7mn+4+16\omega}{32}\right\rfloor,&
	\hbox{$mn\equiv 4\ ({\rm mod }\ 8)$, $\gcd(m,n,2)=1$, \rm{and} $(m,n)\not\in\{(12,3),(3,12)\}$\rm{;}}\\[0.4em]
	
	 7+\frac{\omega}{2}, &\hbox{$(m,n)\in\{(12,3),(3,12)\}$\rm{;}} \\[0.4em]

	\left\lfloor\frac{5mn+4+8\omega}{24}\right\rfloor,
	& \hbox{$mn\equiv 4\ ({\rm mod }\ 8)$ \rm{and} $\gcd(m,n,2)=2$\rm{;}}\\[0.4em]

	\left\lfloor\frac{13mn+40+32\omega}{64}\right\rfloor,&
	\hbox{$mn\equiv 8\ ({\rm mod }\ 16),$  $\gcd(m,n,2)=2$, $(m,n)\not\in\{(2,4),(4,2)\}$, \rm{and} }\\
    &\hbox{$(m,n)\not\in\{(2,12),(4,6),(6,4),(6,12),(12,2),(12,6)\}$ \rm{when} $\lambda_a=3$\rm{;}}\\[0.4em]

    1, &\hbox{$(m,n)\in\{(2,4),(4,2)\}$\rm{;}} \\[0.4em]
	
	\left\lfloor\frac{5mn+8+8\omega}{24}\right\rfloor, &\hbox{$(m,n)\in\{(2,12),(4,6),(6,4),(6,12),(12,2),(12,6)\}$ \rm{and} $\lambda_a=3$\rm{;}} \\[0.4em]
	
	\left\lfloor\frac{13mn+32+32\omega}{64}\right\rfloor,
	&\hbox{$mn\equiv 0\ ({\rm mod }\ 32)$ \rm{and} $\gcd(m,n,2)=2$\rm{,} \rm{or }}\\
	&\hbox{$mn\equiv 16\ ({\rm mod }\ 32)$ \rm{and} $\gcd(m,n,4)=2$\rm{,}}\\
	&\hbox{\rm{except for} $mn\equiv 32\ ({\rm mod }\ 64)$ \rm{and} $\gcd(m,n,8)=4$ \rm{when} $\lambda_a=2$\rm{;}}\\[0.4em]
	
	\frac{13mn-32}{64},
	&\hbox{$mn\equiv 32\ ({\rm mod }\ 64)$, $\gcd(m,n,8)=4$ \rm{and} $\lambda_a=2$\rm{;}}\\[0.4em]

	\left\lfloor\frac{13mn+48+32\omega}{64}\right\rfloor,& \hbox{$mn\equiv 16\ ({\rm mod }\ 32)$ \rm{and} $\gcd(m,n,4)=4$\rm{,}}\\
	& \hbox{\rm{except for} $mn\equiv 144\ ({\rm mod }\ 192)$ \rm{and} $\gcd(m,n,4)=4$ \rm{when} $\lambda_a=2$\rm{;}}\\[0.4em]
	
	\frac{13mn-16}{64},& \hbox{$mn\equiv 144\ ({\rm mod }\ 192)$, $\gcd(m,n,4)=4$ \rm{and} $\lambda_a=2$\rm{.}}\\

	\end{array}
	\right.
	\end{eqnarray}
\end{Theorem}

In Section \ref{add}, we shall establish three recursive constructions for $(m,n,3,\lambda_a,1)$-OOSPCs. Especially, a very efficient doubling construction is presented in Construction \ref{2-mul} to facilitate determining the exact value of $\Theta(m,n,3,\lambda_a,1)$ for $m,n\equiv2\ ({\rm mod }\ 4)$ and $\lambda_a=2,3$. We are to prove the following theorem in Section \ref{4}.
\begin{Theorem}\label{c-main-2} Let $\lambda_a\in\{2,3\}$.
For any $m,n\equiv 2\ ({\rm mod}\ 4)$,
\begin{eqnarray}\nonumber
\Theta(m,n,3,\lambda_a,1)
=\left\lfloor\frac{5mn+4+8\omega}{24}\right\rfloor.
\end{eqnarray}
\end{Theorem}

Finally, in Section \ref{conclusion}, it is conjectured that when $mn\equiv 0\ ({\rm mod }\ 4)$, our bound for $\Theta(m,n,3,\lambda_a,1)$ with $\lambda_a\in\{2,3\}$ shown in Theorem \ref{c-main-1} is tight.

%Note that an optimal $(m,n,k,\lambda_a,\lambda_c)$-OOSPC is equivalent to an optimal $(n,m,k,\lambda_a,\lambda_c)$-OOSPC, hence it holds that $\Theta(m,n,k,\lambda_a,\lambda_c)=\Theta(n,m,k,\lambda_a,\lambda_c)$.

\section{\label{2} Preliminaries}

\subsection{Set-theoretic descriptions}

A convenient way of viewing optical orthogonal signature pattern codes is from a set-theoretic perspective.

Let ${\cal C}$ be an $(m,n,k,\lambda_a,\lambda_c)$-OOSPC. For each $(0,1)$-matrix $M=(a_{ij})\in{\cal C}$,
whose rows are indexed by $\mathbb{Z}_m$ and columns are indexed by $\mathbb{Z}_n$, define
$X_M=\{(i,j)\in \mathbb{Z}_m \times \mathbb{Z}_n:a_{ij}=1\}$. Then, ${\cal{F}}=\{X_M:M\in {\cal C}\}$
is a set-theoretic representation of ${\cal C}$. Conversely, %for any subset $X\subseteq \mathbb{Z}_m\times \mathbb{Z}_n$, construct a $(0,1)$-matrix $X_M=(a_{ij})$ such that $a_{ij}=1$ if and only if $(i,j)\in X$.
let ${\cal F}$ be a set of $k$-subsets of $\mathbb{Z}_m\times \mathbb{Z}_n$. Then ${\cal F}$ is an $(m, n, k, \lambda_a,\lambda_c)$-OOSPC if the following properties are satisfied:
\begin{enumerate}
	\item[{\rm $(1')$}] the auto-correlation property: $|X\cap (X+(s,t))|\leq\lambda_a$ for any $X \in {\cal F}$ and any $(s,t)\in \mathbb{Z}_m\times \mathbb{Z}_n\setminus \{(0,0)\};$
	\item[{\rm $(2')$}] the cross-correlation property: $|X\cap (Y+(s,t))|\leq\lambda_c$ for any
	distinct $X, Y\in {\cal F}$ and any $(s,t)\in \mathbb{Z}_m\times \mathbb{Z}_n,$
\end{enumerate}
where the addition ``+" performs in $\mathbb{Z}_{m}\times \mathbb{Z}_{n}.$
Throughout this paper, we shall use the set-theoretic notation to list codewords of a given OOSPC.

For a given set $\cal F$ of $k$-subsets of $\mathbb{Z}_m\times \mathbb{Z}_n$, it is not convenient to check whether it satisfies the auto- and cross-correlation property according to Conditions $(1')$ and $(2')$. However, when $\lambda_c=1$, a more efficient description can be given by using the difference method. Let $X\in {\cal F}$. Define the list of \emph{differences} of $X$ by
$$\Delta X=\{(x_i,y_i)-(x_j,y_j):0\leq i,j \leq k-1, i\neq j\}$$
as a multiset, and define the \emph{support} of $\Delta X$, denoted by $\textrm{supp}(\Delta X)$,
as the set of underlying elements in $\Delta X$.
%Let $$\Delta{\cal F}=\bigcup\limits_{X\in {\cal F}}\textrm{supp}(\Delta X).$$
Let $\lambda(X)$ denote the maximum multiplicity of elements in the multiset $\Delta X$. Then ${\cal F}$ constitutes an $(m, n, k, \lambda_a,1)$-OOSPC
if the following properties are satisfied:
\begin{enumerate}
\item[{\rm $(1'')$}] the auto-correlation property: $\lambda(X)\leq\lambda_a$ for any $X \in {\cal F}$;
\item[{\rm $(2'')$}] the cross-correlation property: ${\rm{supp}}(\Delta X)\cap {\rm{supp}}(\Delta Y)=\emptyset$ for any distinct $X, Y\in {\cal F}$.
\end{enumerate}

\begin{Example}\label{eg}
We here give an example of a $(6,6,3,\lambda_a,1)$-{\rm OOSPC} with $\lambda_a\in\{2,3\}$ defined on $\mathbb{Z}_6\times \mathbb{Z}_{6}$ as follows:
\begin{center}
\begin{tabular}{llll}
$\lambda_a=2$: & $\{(0,0),(0,3),(3,0)\}$, & $\{(0,0),(0,1),(0,2)\}$, &  $\{(0,0),(1,0),(2,0)\}$, \\ & $\{(0,0),(1,1),(2,2)\}$, &
$\{(0,0),(1,2),(2,1)\}$, & $\{(0,0),(1,3),(3,2)\}$, \\ &  $\{(0,0),(1,4),(3,1)\}$; \\
$\lambda_a=3$: & $\{(0,0),(0,2),(0,4)\}$, & $\{(0,0),(2,0),(4,0)\}$, &  $\{(0,0),(2,2),(4,4)\}$, \\
& $\{(0,0),(2,4),(4,2)\}$, & $\{(0,0),(0,1),(1,0)\}$, & $\{(0,0),(1,1),(2,3)\}$, \\
&  $\{(0,0),(1,3),(3,2)\}$, & $\{(0,0),(1,4),(3,5)\}$, & $\{(0,0),(0,3),(3,0)\}$.
\end{tabular}
\end{center}
\end{Example}

\subsection{Notation and basic propositions}

Throughout this paper, let $A=\mathbb{Z}_m\times \mathbb{Z}_n$. For each   $(x,y)\in A\setminus\{(0,0)\}$, denote by $\pm(x,y)$ the two elements $(x,y)$ and $(-x,-y)$ in $A$.

\begin{Proposition}\label{prop:3-subgroups}
All possible subgroups of order $3$ in $A$ are
\begin{center}
\begin{tabular}{lll}
$\{(0,0),(0,\frac{n}{3}),(0,\frac{2n}{3})\}$, & $\{(0,0)$, $(\frac{m}{3},0),(\frac{2m}{3},0)\}$, \\
$\{(0,0),(\frac{m}{3},\frac{n}{3}),(\frac{2m}{3},\frac{2n}{3})\}$, &
$\{(0,0),(\frac{m}{3},\frac{2n}{3}),(\frac{2m}{3},\frac{n}{3})\}$.
 \end{tabular}
 \end{center}
\end{Proposition}

\begin{Proposition}\label{prop:4-subgroups}
All possible cyclic subgroups of order $4$ in $A$ are
\begin{center}
\begin{tabular}{lll}
$\{(0,0),\pm(\frac{m}{4},0),(\frac{m}{2},0)\}$, &
 $\{(0,0),\pm(\frac{m}{4},\frac{n}{4}), (\frac{m}{2},\frac{n}{2})\}$, &
$\{(0,0), \pm(\frac{m}{4},\frac{n}{2}),(\frac{m}{2},0)\}$, \\
$\{(0,0), \pm(\frac{m}{4},\frac{3n}{4}), (\frac{m}{2},\frac{n}{2})\}$, &
$\{(0,0),\pm(0,\frac{n}{4}),(0,\frac{n}{2})\}$, &
$\{(0,0), \pm(\frac{m}{2},\frac{n}{4}), (0,\frac{n}{2})\}$.
 \end{tabular}
 \end{center}
\end{Proposition}

\begin{Proposition}\label{prop:4-subgroups-klein}
The unique possible subgroup of order $4$ isomorphic to $\mathbb Z_2\times \mathbb Z_2$ in $A$ is
$\{(0,0),(\frac{m}{2},0),(0,\frac{n}{2}),(\frac{m}{2},\frac{n}{2})\}$.
\end{Proposition}

%We explain some notations that will be used in the next section.

Let $(G,+)$ be an abelian group with the identity $0$. Let $X\subseteq G$. The $G$-\emph{orbit} of $X$ is the set ${\rm Orb}_G(X)=\{X+g:g\in G\}$, where $X+g=\{x+g: x \in X\}$.
%For $\alpha\in G$, the \emph{order} of $\alpha$, briefly \emph{${\rm ord}_G(\alpha)$}, is the smallest positive integer $i$ such that $i\alpha=0$.
For any positive integer $i$, let
%The \emph{stabilizer} of $X$ \emph{under} $G$ is a subgroup of $G$ consisting of all elements $\beta\in G$ such that $X+\beta=X$, where $X+\beta=\{\alpha+\beta:\alpha \in X\}$.
%$\Omega_G(j)$ denote the set of all the elements in $A$ such that $i$ can be divided by their orders, i.e.,
%$$\Omega_G(j)=\{\alpha\in G:j \alpha=0\}=\{\alpha \in G:{\rm ord}_G(\alpha)\mid j\}.$$
$$\Omega_G(i)=\{\alpha\in G:i \alpha=0\}.$$

\begin{Proposition}\label{prop:omiga2}
$\Omega_A(2)=
\left\{(0,0),\left(0,\frac{n}{2}\right),\left(\frac{m}{2},0\right),
\left(\frac{m}{2},\frac{n}{2}\right)\right\}$.
\end{Proposition}

\begin{Remark}\label{rmk-omega}
In Proposition $\ref{prop:omiga2}$, the notation $\Omega_A(2)$ should be understood as follows
$$\Omega_A(2)=\left\{
\begin{array}{ll}
\left\{(0,0)\right\},&   {\rm if}\  m,n\equiv 1\ ({\rm mod }\ 2);\\
\left\{(0,0),\left(\frac{m}{2},0\right)\right\},&  {\rm if}\  m\equiv 0\ ({\rm mod }\ 2)\ \hbox{and}\ n\equiv 1\ ({\rm mod }\ 2);\\
\left\{(0,0),\left(0,\frac{n}{2}\right)\right\},&  {\rm if}\  m\equiv 1\ ({\rm mod }\ 2)\ \hbox{and}\ n\equiv 0\ ({\rm mod }\ 2);\\
\left\{(0,0),\left(0,\frac{n}{2}\right),\left(\frac{m}{2},0\right),\left(\frac{m}{2},\frac{n}{2}\right)\right\},&
{\rm if}\  m,n\equiv 0\ ({\rm mod }\ 2).
\end{array}\right.$$
In what follows, we always use a similar method to denote sets. The reader can judge it according to the context.
\end{Remark}

\begin{Proposition}\label{prop:omiga3}
\begin{itemize}
\item[$(1)$] $\Omega_A(3)=\left\{(0,0),\pm\left(0,\frac{n}{3}\right),\pm\left(\frac{m}{3},0\right),\pm\left(\frac{m}{3},\frac{n}{3}\right),
\pm\left(\frac{m}{3},\frac{2n}{3}\right)\right\}.
$
\item[$(2)$] $\Omega_A(4)=\Omega_A(2)\cup \left\{\pm\left(0,\frac{n}{4}\right),\pm\left(\frac{m}{4},0\right),\pm\left(\frac{m}{4},\frac{n}{4}\right),
\pm\left(\frac{m}{4},\frac{n}{2}\right),\pm\left(\frac{m}{4},\frac{3n}{4}\right),
\pm\left(\frac{m}{2},\frac{n}{4}\right)\right\}.
$
\end{itemize}
\end{Proposition}

\section{\label{331} Upper bound on the size of $(m, n, 3, \lambda_a,1)$-OOSPCs}

In this section, we shall estimate the upper bound of $\Theta(m,n,3,\lambda_a,1)$ with $\lambda_a\in\{2,3\}$ for any positive integers $m$ and $n$. Without loss of generality assume that each codeword in an OOSPC contains the element $(0,0)$. Let
\begin{itemize}
\item[] ${\cal{T}}_1=\left\{{\rm Orb}_A(\{(0,0),\alpha,2\alpha\}):\alpha\in  A\setminus (\Omega_A(3)\cup\Omega_A(4))\right\}$, and
\item[] ${\cal{T}}_2=\left\{{\rm Orb}_A(\{(0,0),\alpha,\beta\}):\alpha\in  A\setminus \Omega_A(4)~{\rm and}~\beta\in \Omega_A(2)\setminus\{(0,0)\}\right\}.
$
\end{itemize}

\noindent The following two lemmas play key roles to derive our bound.

\begin{Lemma} {\rm \cite[Lemma 4.1]{sawa}}
\label{0-1} Let $X=\{(0,0),\alpha,\beta\}$ be a $3$-subset of $A$.
Then
\begin{eqnarray}
|{\rm{supp}}(\Delta X)|
=\left\{
\begin{array}{ll}
2,\ \hbox{{\rm if} $<\alpha,\beta>\cong \mathbb{Z}_3${\rm ;}}\\
3,\ \hbox{{\rm if} $<\alpha,\beta>\cong \mathbb{Z}_4$ or $\mathbb{Z}_2\times \mathbb{Z}_2${\rm ;}}\\
4,\ \hbox{{\rm if} ${\rm Orb}_A(X)\in{\cal T}_1${\rm ;}}\\
5,\ \hbox{{\rm if} ${\rm Orb}_A(X)\in{\cal T}_2${\rm ;}}\\
6,\ \hbox{{\rm otherwise,}}\\
\end{array}
\right.\nonumber
\end{eqnarray}
where $<\alpha,\beta>$ denotes the additive subgroup of $A$ generated by $\alpha$ and $\beta$.
\end{Lemma}

\begin{Lemma}{\rm \cite[Lemma 4.2]{sawa}}
\label{0-2} Let $X$ be a $3$-subset of $A$.
	Then
	\begin{eqnarray}
	\lambda(X)=\max \limits_{(s,t)\in A\setminus \{(0,0)\}} |X\cap (X+(s,t))|
	=\left\{
	\begin{array}{ll}
	3,\ \hbox{{\rm if} $|{\rm{supp}}(\Delta X)|=2${\rm ;}}\\
	2,\ \hbox{{\rm if} $|{\rm{supp}}(\Delta X)|=3,4,5${\rm ;}}\\
	1,\ \hbox{{\rm if} $|{\rm{supp}}(\Delta X)|=6${\rm .}}\\
	\end{array}
	\right.\nonumber
	\end{eqnarray}
\end{Lemma}

\subsection{General upper bound}

Let $\lambda_a\in\{2,3\}$. For any codeword $X$ of an $(m, n, 3,\lambda_a,1)$-OOSPC $\cal F$, if $|{\rm{supp}}(\Delta X)|=i$, then $X$ is said to be of \emph{Type $i$}. By Lemma \ref{0-1}, $i\in\{2,3,4,5,6\}$. Let $N_{i}$ denote the number of codewords in $\cal F$ of Type $i$.
The cross-correlation property $(2'')$ implies that $\Delta{\cal F}=\bigcup_{X\in \cal{F}}{\rm{supp}}(\Delta X)$ covers each nonzero element of $A$ at most once. Thus we have
\begin{eqnarray}\label{mn-1}
2N_2+3N_3+4N_4+5N_5+6N_6\leq mn-1.
\end{eqnarray}

\begin{Lemma}\label{1-1}
\begin{eqnarray}\label{N2}
	N_2\leq \omega.
	\end{eqnarray}
\end{Lemma}

\proof Let $X=\{(0,0),\alpha,\beta\}$ be a codeword satisfying $|{\rm{supp}}(\Delta X)|=2$. By Lemma \ref{0-2}, $\lambda(X)=3$. Note that the auto-correlation property $(1'')$ requires $\lambda(X)\leq\lambda_a$ for any $X \in {\cal F}$. So if $\lambda_a=2$, then $N_2=0$. If $\lambda_a=3$, by Lemma \ref{0-1}, $<\alpha,\beta>$ forms an additive subgroup of order $3$ in $A$. Since all possible subgroups of order $3$ in $A$ are $\{(0,0),(0,\frac{n}{3}),(0,\frac{2n}{3})\}$, $\{(0,0),(\frac{m}{3},0),(\frac{2m}{3},0\}$, $\{(0,0),(\frac{m}{3},\frac{n}{3}),(\frac{2m}{3},\frac{2n}{3})\}$ and  $\{(0,0),(\frac{m}{3},\frac{2n}{3}),(\frac{2m}{3},\frac{n}{3})\}$, the value of $N_2$ depends on whether $m$ and $n$ could be divided by $3$. \qed

For each codeword $X=\{(0,0),\alpha,\beta\}$ of Type $3$, if $<\alpha,\beta>\cong \mathbb{Z}_4$, then by Proposition \ref{prop:4-subgroups}, w.l.o.g., $X$ is one of the following forms:
\begin{eqnarray}\label{sub1}
	\begin{tabular}{lll}
$\{(0,0),(\frac{m}{4},0),(\frac{m}{2},0)\}$, & $\{(0,0),(\frac{m}{4},\frac{n}{4}),(\frac{m}{2},\frac{n}{2})\}$, &
$\{(0,0),(\frac{m}{4},\frac{n}{2}),(\frac{m}{2},0)\}$; \\
$\{(0,0),(\frac{m}{4},\frac{3n}{4}),(\frac{m}{2},\frac{n}{2})\}$, &
$\{(0,0),(0,\frac{n}{4}),(0,\frac{n}{2})\}$, &
$\{(0,0),(\frac{m}{2},\frac{n}{4}),(0,\frac{n}{2})\}$.
\end{tabular}
\end{eqnarray}
If $<\alpha,\beta>\cong \mathbb{Z}_2\times \mathbb Z_2$, then by Proposition \ref{prop:4-subgroups-klein}, w.l.o.g., $X$ is of the form
\begin{eqnarray}\label{sub2}
\{(0,0),(\frac{m}{2},0),(0,\frac{n}{2})\}.
\end{eqnarray}
Therefore, all codewords of Type $3$, which are in the form of $\{(0,0),\alpha=(a,b),\beta\}$, can be divided into the following three types:
\begin{itemize}
\item[]\textbf{Type 3.1}: $\alpha,\beta\in\Omega_A(2)\setminus\{(0,0)\}$;
\item[]\textbf{Type 3.2}: $\alpha\in \Omega_{A}(4)\setminus\{(0,0)\}$ and $a\in \Omega_{\mathbb{Z}_m}(4)\setminus\Omega_{\mathbb{Z}_m}(2)$;
\item[]\textbf{Type 3.3}: $\alpha\in \Omega_{A}(4)\setminus\{(0,0)\}$, $a\in\Omega_{\mathbb{Z}_m}(2)$ and $b\in\Omega_{\mathbb{Z}_n}(4)\setminus \Omega_{\mathbb{Z}_n}(2)$.
\end{itemize}

\noindent Let $N_3^{(1)}$, $N_3^{(2)}$ and $N_3^{(3)}$ denote the number of codewords in $\cal{F}$ of Types $3.1$, $3.2$ and $3.3$, respectively.
Then, $$N_3=N_3^{(1)}+N_3^{(2)}+N_3^{(3)}.$$

\begin{Remark}\label{for type3.1.2.3}
By \eqref{sub1} and \eqref{sub2}, one can check the following facts.
\begin{enumerate}
\item[$(1)$] For any codeword $X$ of Type $3.1$, ${\rm supp}(\Delta X)$ is of the form  $\{(\frac{m}{2},0), (0,\frac{n}{2}), (\frac{m}{2},\frac{n}{2})\}.$
\item[$(2)$] For any codeword $X$ of Type $3.2$, ${\rm supp}(\Delta X)$ is one of the forms:
$\{\pm(\frac{m}{4},0), (\frac{m}{2},0)\}$, $\{\pm(\frac{m}{4}$, $\frac{n}{4}),(\frac{m}{2},\frac{n}{2})\}$, $\{\pm(\frac{m}{4},\frac{n}{2}), (\frac{m}{2},0)\}$ and $\{\pm(\frac{m}{4},\frac{3n}{4}),(\frac{m}{2},\frac{n}{2})\}.
$
\item[$(3)$] For any codeword $X$ of Type $3.3$, ${\rm supp}(\Delta X)$ is one of the forms:
$\{\pm(0,\frac{n}{4}), (0,\frac{n}{2})\}$ and $\{\pm(\frac{m}{2},\frac{n}{4}), (0,\frac{n}{2})\}$.
\end{enumerate}
\end{Remark}

\begin{Lemma}\label{1-4}\begin{eqnarray}\label{N31}
	N_3^{(1)}\leq \left\{
	\begin{array}{ll}
	1,\ \hbox{ \rm{if} $\gcd(m,n,2)=2$\rm{;}}\\
	0,\ \hbox{ \rm{if} $\gcd(m,n,2)=1$\rm{.}}\\
	\end{array}
	\right.
	\end{eqnarray}
\end{Lemma}
\proof For each codeword $X$ of Type $3.1$, $|{\rm{supp}}(\Delta X)\cap\Omega_A(2)|=3$ by Remark \ref{for type3.1.2.3}(1). By Remark \ref{rmk-omega}, we have $$3N_3^{(1)}\leq |\Omega_A(2)\setminus\{(0,0)\}|=\left\{
	\begin{array}{ll}
	3,\ \hbox{ \rm{if} $\gcd(m,n,2)=2$\rm{;}}\\
	0,\ \hbox{ \rm{if} $\gcd(m,n,2)=1$\rm{.}}\\
	\end{array}
	\right.$$ \qed

\begin{Lemma}\label{1-5}
Let
\begin{eqnarray}\nonumber
	\varepsilon=\left\{
	\begin{array}{ll}
	2,\ \hbox{ \rm{if} $m,n\equiv 0\ ({\rm mod }\ 4)$\rm{;}}\\
	1,\ \hbox{ \rm{if} $m\equiv 0\ ({\rm mod }\ 4)$ \rm{and} $n\not\equiv 0\ ({\rm mod }\ 4)$\rm{;}}\\
	0,\ \hbox{ \rm{if} $m\not\equiv 0\ ({\rm mod }\ 4)$\rm{.}}\\
	\end{array}
	\right.
\end{eqnarray}
Then
\begin{flalign}
\quad\label{N32} (1)\ \ N_3^{(2)}&\leq \varepsilon; \\
\quad\label{N33} (2)\ \ N_3^{(3)}&\leq \left\{
	\begin{array}{ll}
	1,\ \hbox{ \rm{if} $n\equiv 0\ ({\rm mod }\ 4)$\rm{;}}\\
	0,\ \hbox{ \rm{if} $n\not\equiv 0\ ({\rm mod }\ 4)$\rm{.}}\\
	\end{array}
	\right.&
\end{flalign}
\end{Lemma}

\proof Each codeword of Type $3$ is one of the forms shown in \eqref{sub1} and \eqref{sub2}.

$(1)$ When $m\not\equiv 0\ ({\rm mod }\ 4)$, there is no codeword of Type $3.2$, and so $N_3^{(2)}=0$.
When $m\equiv 0\ ({\rm mod }\ 4)$ and $n\equiv 1\ ({\rm mod }\ 2)$, there is at most one codeword of Type $3.2$, which is $\{(0,0),(\frac{m}{4},0),(\frac{m}{2},0)\}$, and so $N_3^{(2)}\leq1$. When $m\equiv 0\ ({\rm mod }\ 4)$ and $n\equiv 2\ ({\rm mod }\ 4)$, there are at most two codewords of Type $3.2$, which are $\{(0,0),(\frac{m}{4},0),(\frac{m}{2},0)\}$ and $\{(0,0),(\frac{m}{4},\frac{n}{2}),(\frac{m}{2},0)\}$. Here $(\frac{m}{2},0)$ is shared as a difference. Hence, $N_3^{(2)}\leq 1$.
When $m,n\equiv 0\ ({\rm mod }\ 4)$, there are at most four codewords of Type $3.2$, which are $\{(0,0),(\frac{m}{4},0),(\frac{m}{2},0)\},\{(0,0),(\frac{m}{4}$, $\frac{n}{4}),(\frac{m}{2},\frac{n}{2})\}$,
$\{(0,0),(\frac{m}{4},\frac{n}{2}),(\frac{m}{2},0)\}$ and
$\{(0,0),(\frac{m}{4},\frac{3n}{4}),(\frac{m}{2},\frac{n}{2})\}$. Since $\lambda_c=1$, it is readily checked that $N_3^{(2)}\leq2$.

$(2)$ When $n\not\equiv 0\ ({\rm mod }\ 4)$, there is no codeword of Type $3.3$, and so $N_3^{(3)}=0$.
When $n\equiv 0\ ({\rm mod }\ 4)$, there are at most two codewords of Type $3.3$, which are
$\{(0,0),(0,\frac{n}{4}),(0,\frac{n}{2})\}$ and $\{(0,0),(\frac{m}{2},\frac{n}{4}),(0,\frac{n}{2})\}$. Here $(0,\frac{n}{2})$ is shared as a difference. Hence, $N_3^{(3)}\leq 1$. \qed

\begin{Lemma}\label{odd-m}
$\Theta(m,n,3,\lambda_a,1)\leq \left\lfloor\frac{ mn+2\omega}{4}\right\rfloor$ for any $m\equiv 1\ ({\rm mod }\ 2)$ and $\lambda_a\in\{2,3\}$.
\end{Lemma}

\proof For $m\equiv 1\ ({\rm mod }\ 2)$, by (\ref{N2}), (\ref{N31}), (\ref{N32}) and (\ref{N33}), we have $N_2\leq\omega$, $N_3^{(1)}=N_3^{(2)}=0$ and $N_3^{(3)}\leq1$. Then by (\ref{mn-1}), $4(N_2+N_3+N_4+N_5+N_6)\leq mn-1+2N_2+N_3-N_5-2N_6\leq mn+2\omega$. Thus $\Theta(m,n,3,\lambda_a,1)\leq \left\lfloor\frac{ mn+2\omega}{4}\right\rfloor$. \qed

\begin{Lemma}\label{1-7}
Let
\begin{eqnarray}\nonumber
	\rho=\left\{
	\begin{array}{ll}
	3,\ \hbox{ \rm{if} $\gcd(m,n,2)=2$\rm{;}}\\
	1,\ \hbox{ \rm{if} $2\mid mn$ \rm{and} $\gcd(m,n,2)=1$\rm{;}}\\
	0,\ \hbox{ \rm{if} $2\nmid mn$\rm{.}}\\
	\end{array}
	\right.
	\end{eqnarray}
Then
\begin{eqnarray}\label{N3andN5}
	3N_3^{(1)}+N_3^{(2)}+N_3^{(3)}+N_5\leq \rho.
	\end{eqnarray}
\end{Lemma}
\proof By Proposition \ref{prop:omiga2}, $\Omega_A(2)=\{(0,0),(0,\frac{n}{2}),(\frac{m}{2},0),(\frac{m}{2},\frac{n}{2})\}$. If $X$ is a codeword of Type $3.1$, then $|{\rm{supp}}(\Delta X)\cap\Omega_A(2)|=3$ by Remark \ref{for type3.1.2.3}(1). If $X$ is a codeword of Types $3.2$ or $3.3$ or Type $5$, then $|{\rm{supp}}(\Delta X)\cap\Omega_A(2)|=1$ by Remark \ref{for type3.1.2.3}(2)(3) and Lemma \ref{0-1}. Hence, $3N_3^{(1)}+N_3^{(2)}+N_3^{(3)}+N_5\leq|\Omega_A(2)\setminus\{(0,0)\}|=\rho$. \qed

Without loss of generality, each codeword of Type $4$ is of the form $\{(0,0),(a,b),(2a,2b)\}$, where $(a,b)\in A\setminus(\Omega_A(3)\cup\Omega_A(4))$. We divide the codewords of Type $4$ into the following two types according to the parity of $a$:
\begin{itemize}
\item[]\textbf{Type 4.1}: $a\equiv 1\ ({\rm mod }\ 2)$;
\item[]\textbf{Type 4.2}: $a\equiv 0\ ({\rm mod }\ 2)$.
\end{itemize}
\noindent Let $N_4^{(1)}$ and $N_4^{(2)}$ denote the number of codewords in $\cal{F}$ of Types $4.1$ and $4.2$, respectively.
Then, $$N_4=N_4^{(1)}+N_4^{(2)}.$$

Throughout this paper, we always set
\begin{center}
$\begin{array}{ll}
  A_{e\cdot}=\left\{(x,y)\in A:\ x\equiv 0\ ({\rm mod }\ 2)\right\}, & A_{eo}=\left\{(x,y)\in A:\ x\equiv 0\ ({\rm mod }\ 2),\ y\equiv 1\ ({\rm mod }\ 2)\right\}, \\[0.5em]
  A_{s\cdot}=\left\{(x,y)\in A:\ x\equiv 2\ ({\rm mod }\ 4)\right\}, & A_{oe}=\left\{(x,y)\in A:\ x\equiv 1\ ({\rm mod }\ 2),\ y\equiv 0\ ({\rm mod }\ 2)\right\}, \\[0.5em]
  A_{d\cdot}=\left\{(x,y)\in A:\ x\equiv 0\ ({\rm mod }\ 4)\right\}, & A_{se}=\left\{(x,y)\in A:\ x\equiv 2\ ({\rm mod }\ 4),\ y\equiv 0\ ({\rm mod }\ 2)\right\}, \\[0.5em]
 A_{\cdot o}=\left\{(x,y)\in A:\ y\equiv 1\ ({\rm mod }\ 2)\right\}, & A_{de}=\left\{(x,y)\in A:\ x\equiv 0\ ({\rm mod }\ 4),\ y\equiv 0\ ({\rm mod }\ 2)\right\}, \\[0.5em]
A_{ee}=\left\{(x,y)\in A:\ x,y\equiv 0\ ({\rm mod }\ 2)\right\}, &A_{ds}=\left\{(x,y)\in A:\ x\equiv 0\ ({\rm mod }\ 4),\ y\equiv 2\ ({\rm mod }\ 4)\right\}, \\[0.5em]
A_{oo}=\left\{(x,y)\in A:\ x,y\equiv 1\ ({\rm mod }\ 2)\right\}.\\[0.5em]
\end{array}$
\end{center}

\begin{Lemma}\label{1-3}
\begin{eqnarray}\label{N2andN4}
2N_2+2N_4\leq \left\{
\begin{array}{ll}
\frac{mn}{4}-1,& \hbox{ \rm{if} $m,n\equiv 2\ ({\rm mod }\ 4)$\rm{;}}\\[0.2em]
\frac{mn}{4}-2,& \hbox{ \rm{if} $m\equiv 2\ ({\rm mod }\ 4)$ {\rm and} $n\equiv 0\ ({\rm mod }\ 4)$,}\\
& \hbox{ {\rm or} $m\equiv 0\ ({\rm mod }\ 4)$ {\rm and} $n\equiv 2\ ({\rm mod }\ 4)$\rm{;}}\\[0.2em]
\frac{mn}{4}-4,& \hbox{ \rm{if} $m,n\equiv 0\ ({\rm mod }\ 4)$\rm{.}}
\end{array}
\right.
\end{eqnarray}
\end{Lemma}
\proof Due to $|A_{ee}|=\frac{m}{2}\times\frac{n}{2}=\frac{mn}{4}$ and
\begin{eqnarray}\nonumber
A_{ee}\cap\Omega_A(2)=\left\{
\begin{array}{ll}
\{(0,0)\},& \hbox{ \rm{if} $m,n\equiv 2\ ({\rm mod }\ 4)$\rm{;}}\\[0.2em]
\{(0,0),(0,\frac{n}{2})\},& \hbox{ \rm{if} $m\equiv 2\ ({\rm mod }\ 4)$ {\rm and} $n\equiv 0\ ({\rm mod }\ 4)$\rm{;}}\\[0.2em]
\{(0,0),(\frac{m}{2},0)\},& \hbox{ \rm{if} $m\equiv 0\ ({\rm mod }\ 4)$ {\rm and} $n\equiv 2\ ({\rm mod }\ 4)$\rm{;}}\\[0.2em]
\{(0,0),(0,\frac{n}{2}),(\frac{m}{2},0),(\frac{m}{2},\frac{n}{2})\},& \hbox{ \rm{if} $m,n\equiv 0\ ({\rm mod }\ 4)$,}\\
\end{array}
\right.
\end{eqnarray}
we have
\begin{eqnarray}\nonumber
|A_{ee}\setminus\Omega_A(2)|=|A_{ee}|-|A_{ee}\cap\Omega_A(2)|=\left\{
\begin{array}{ll}
\frac{mn}{4}-1,& \hbox{ \rm{if} $m,n\equiv 2\ ({\rm mod }\ 4)$\rm{;}}\\[0.2em]
\frac{mn}{4}-2,& \hbox{ \rm{if} $m\equiv 2\ ({\rm mod }\ 4)$ {\rm and} $n\equiv 0\ ({\rm mod }\ 4)$,}\\ & \hbox{ {\rm or} $m\equiv 0\ ({\rm mod }\ 4)$ {\rm and} $n\equiv 2\ ({\rm mod }\ 4)$\rm{;}}\\[0.2em]
\frac{mn}{4}-4,& \hbox{ \rm{if} $m,n\equiv 0\ ({\rm mod }\ 4)$\rm{.}}\\
\end{array}
\right.
\end{eqnarray}
Since $m$ and $n$ are even, it is readily checked that for any codeword $X$ of Type $2$, ${\rm{supp}}(\Delta X)\subseteq A_{ee}\setminus\Omega_A(2)$. For each codeword $X=\{(0,0),(a,b),(2a,2b)\}$ of Type $4$, ${\rm{supp}}(\Delta X)$ contributes at least two differences, $\pm(2a,2b)$, in $A_{ee}\setminus\Omega_A(2)$. Thus $2N_2+2N_4\leq|A_{ee}\setminus\Omega_A(2)|$.
\qed

\begin{Lemma}\label{1-8}
For any $m\equiv 2\ ({\rm mod }\ 4)$ and $\lambda_a\in\{2,3\}$,
\begin{eqnarray}\label{e-8}\nonumber
\Theta(m,n,3,\lambda_a,1)\leq\left\{
\begin{array}{ll}
\left\lfloor\frac{5mn+4+8\omega}{24}\right\rfloor, & \hbox{ \rm{if} $n\equiv 2\ ({\rm mod }\ 4)$\rm{;}}\\[0.4em]
\left\lfloor\frac{5mn+8+8\omega}{24}\right\rfloor, & \hbox{ \rm{if} $n\equiv 0\ ({\rm mod }\ 4)$\rm{.}}\\
\end{array}
\right.
\end{eqnarray}
\end{Lemma}
\proof For $m\equiv 2\ ({\rm mod }\ 4)$, by \eqref{N32}, $N_3^{(2)}=0$ and
hence $N_3=N_3^{(1)}+N_3^{(3)}$. We rewrite \eqref{mn-1} and \eqref{N3andN5} as follows:
\begin{align}
\label{e-11} 2N_2+3N_3^{(1)}+3N_3^{(3)}+4N_4+5N_5+6N_6&\leq mn-1, \\
\label{e-12} 3N_3^{(1)}+N_3^{(3)}+N_5&\leq 3.
\end{align}
By $2\times$(\ref{N2})+$2\times$(\ref{N33})+(\ref{N2andN4})+(\ref{e-11})+(\ref{e-12}), we have
\begin{eqnarray}\nonumber
6(N_2+N_3+N_4+N_5+N_6)\leq\left\{
\begin{array}{ll}
\frac{5mn}{4}+2\omega+1, & \hbox{ \rm{if} $n\equiv 2\ ({\rm mod }\ 4)$\rm{;}}\\[0.2em]
\frac{5mn}{4}+2\omega+2, & \hbox{ \rm{if} $n\equiv 0\ ({\rm mod }\ 4)$\rm{.}}\\
\end{array}
\right.
\end{eqnarray}
Since $\Theta(m,n,3,\lambda_a,1)=N_2+N_3+N_4+N_5+N_6$, the conclusion follows. \qed

\begin{Remark}\label{1-r-1}
Examining the proof of Lemma $\ref{1-8}$, we have that if
\begin{eqnarray}\label{e-8}\nonumber
\Theta(m,n,3,\lambda_a,1)=\left\{
\begin{array}{ll}
\left\lfloor\frac{5mn+4+8\omega}{24}\right\rfloor, & \hbox{ \rm{if} $n\equiv 2\ ({\rm mod }\ 4)$\rm{;}}\\[0.4em]
\left\lfloor\frac{5mn+8+8\omega}{24}\right\rfloor, & \hbox{ \rm{if} $n\equiv 0\ ({\rm mod }\ 4)$\rm{,}}\\
\end{array}
\right.
\end{eqnarray}
where $m\equiv 2\ ({\rm mod }\ 4)$ and $\lambda_a\in\{2,3\}$,
then the equalities must hold in $(\ref{N2})$, $(\ref{N33})$, $(\ref{N2andN4})$, $(\ref{e-11})$ and $(\ref{e-12})$.
\end{Remark}

\begin{Lemma}\label{1-2}
Let
\begin{eqnarray}\nonumber
	\eta=\left\{
	\begin{array}{ll}
	\frac{mn}{8},& \hbox{ \rm{if }$m\equiv 0\ ({\rm mod }\ 8)$ \rm{and} $n\equiv 0\ ({\rm mod }\ 2)$\rm{;}}\\[0.2em]
	\frac{mn}{4},& \hbox{  \rm{if} $m\equiv 0\ ({\rm mod }\ 8)$ \rm{and} $n\equiv 1\ ({\rm mod }\ 2)$\rm{;}}\\[0.2em]
	\frac{mn}{8}-2,& \hbox{ \rm{if}  $m\equiv 4\ ({\rm mod }\ 8)$ \rm{and} $n\equiv 0\ ({\rm mod }\ 4)$\rm{;}}\\[0.2em]
	\frac{mn}{8}-1,& \hbox{  \rm{if} $m\equiv 4\ ({\rm mod }\ 8)$ \rm{and} $n\equiv 2\ ({\rm mod }\ 4)$\rm{;}}\\[0.2em]
	\frac{mn}{4}-1,& \hbox{  \rm{if} $m\equiv 4\ ({\rm mod }\ 8)$ \rm{and} $n\equiv 1\ ({\rm mod }\ 2)$.}\\
	\end{array}
	\right.
	\end{eqnarray}
Then
\begin{eqnarray}\label{N41}
	2N_4^{(1)}\leq \eta.
	\end{eqnarray}
\end{Lemma}
\proof \noindent\textbf{Case 1}: $m\equiv 0\ ({\rm mod }\ 4)$ and $n\equiv 1\ ({\rm mod }\ 2)$.

Due to $|A_{s\cdot}|=\frac{m}{4}\times n=\frac{mn}{4}$ and
\begin{eqnarray}
	A_{s\cdot}\cap \Omega_A(2)=\left\{
	\begin{array}{ll}
	\emptyset,& \hbox{  \rm{if} $m\equiv 0\ ({\rm mod }\ 8)$ \rm{and} $n\equiv 1\ ({\rm mod }\ 2)$\rm{;}}\\[0.2em]
 \{(\frac{m}{2},0)\},& \hbox{  \rm{if} $m\equiv 4\ ({\rm mod }\ 8)$ \rm{and} $n\equiv 1\ ({\rm mod }\ 2)$,}\\
	\end{array}
	\right.\nonumber
	\end{eqnarray}
we have
\begin{eqnarray}
	|A_{s\cdot}\setminus \Omega_A(2)|=|A_{s\cdot}|-|A_{s\cdot}\cap \Omega_A(2)|=\left\{
	\begin{array}{ll}
	\frac{mn}{4},& \hbox{  \rm{if} $m\equiv 0\ ({\rm mod }\ 8)$ \rm{and} $n\equiv 1\ ({\rm mod }\ 2)$\rm{;}}\\[0.2em]
	\frac{mn}{4}-1,& \hbox{  \rm{if} $m\equiv 4\ ({\rm mod }\ 8)$ \rm{and} $n\equiv 1\ ({\rm mod }\ 2)$.}\\
	\end{array}
	\right.\nonumber
	\end{eqnarray}
Each codeword $\{(0,0),(a,b),(2a,2b)\}$ of Type $4.1$ contributes exactly two differences, $\pm(2a,2b)$, in $A_{s\cdot}\setminus \Omega_A(2)$.
Hence, $2N_4^{(1)}\leq |A_{s\cdot}\setminus \Omega_A(2)|$.

\noindent\textbf{Case 2}: $m\equiv 0\ ({\rm mod }\ 4)$ and $n\equiv 0\ ({\rm mod }\ 2)$.

Due to $|A_{se}|=\frac{m}{4}\times \frac{n}{2}=\frac{mn}{8}$ and
\begin{eqnarray}
	A_{se}\cap \Omega_A(2)=\left\{
	\begin{array}{ll}
	\emptyset,& \hbox{  \rm{if} $m\equiv 0\ ({\rm mod }\ 8)$ \rm{and} $n\equiv 0\ ({\rm mod }\ 2)$\rm{;}}\\[0.2em]
    \{(\frac{m}{2},0),(\frac{m}{2},\frac{n}{2})\},& \hbox{  \rm{if} $m\equiv 4\ ({\rm mod }\ 8)$ \rm{and} $n\equiv 0\ ({\rm mod }\ 4)$\rm{;}}\\[0.2em]
    \{(\frac{m}{2},0)\},& \hbox{  \rm{if} $m\equiv 4\ ({\rm mod }\ 8)$ \rm{and} $n\equiv 2\ ({\rm mod }\ 4)$,}\\
	\end{array}
	\right.\nonumber
	\end{eqnarray}
we have
\begin{eqnarray}
	|A_{se}\setminus \Omega_A(2)|=|A_{se}|-|A_{se}\cap \Omega_A(2)|=\left\{
	\begin{array}{ll}
	\frac{mn}{8},& \hbox{ \rm{if }$m\equiv 0\ ({\rm mod }\ 8)$ \rm{and} $n\equiv 0\ ({\rm mod }\ 2)$\rm{;}}\\[0.2em]
	\frac{mn}{8}-2,& \hbox{ \rm{if}  $m\equiv 4\ ({\rm mod }\ 8)$ \rm{and} $n\equiv 0\ ({\rm mod }\ 4)$\rm{;}}\\[0.2em]
	\frac{mn}{8}-1,& \hbox{  \rm{if} $m\equiv 4\ ({\rm mod }\ 8)$ \rm{and} $n\equiv 2\ ({\rm mod }\ 4)$\rm{.}}\\
	\end{array}
	\right.\nonumber
	\end{eqnarray}
Each codeword $\{(0,0),(a,b),(2a,2b)\}$ of Type $4.1$ contributes exactly two differences, $\pm(2a,2b)$, in $A_{se}\setminus \Omega_A(2)$.
Hence, $2N_4^{(1)}\leq |A_{se}\setminus \Omega_A(2)|$. \qed

\begin{Lemma}\label{1-9}
Let \begin{eqnarray}
	\gamma=\left\{
	\begin{array}{ll}
	1,\ \hbox{ \rm{if }$m\equiv 4\ ({\rm mod }\ 8)$\rm{;}}\\
	3,\ \hbox{ \rm{if }$m\equiv 0\ ({\rm mod }\ 8)$\rm{.}}\\
	\end{array}
	\right.\nonumber
	\end{eqnarray}
Then
\begin{eqnarray}\label{mn/2-1}
2N_2+3N_3^{(1)}+\gamma N_3^{(2)}+3N_3^{(3)}+2N_4^{(1)}+4N_4^{(2)}+N_5+2N_6\leq \frac{mn}{2}-1.
\end{eqnarray}

\end{Lemma}
\proof For $m\equiv 0\ ({\rm mod }\ 4)$, $|A_{e\cdot}|=\frac{m}{2}\times n=\frac{mn}{2}$. Take any codeword $X$ of an $(m, n, 3, \lambda_a,1)$-OOSPC, ${\cal F}$, with $\lambda_a\in\{2,3\}$. It is readily checked that
\begin{eqnarray}
	|{\rm{supp}}(\Delta X)\cap A_{e\cdot}|\geq\left\{
	\begin{array}{ll}
	1,\ \hbox{ \rm{if} $X$ \rm{is of Type 5}\rm{;}}\\
	2,\ \hbox{ \rm{if} $X$ \rm{is of Type 2 or 4.1 or 6}\rm{;}}\\
	3,\ \hbox{  \rm{if} $X$ \rm{is of Type 3.1 or 3.3}\rm{;}}\\
	4,\ \hbox{  \rm{if} $X$ \rm{is of Type 4.2}\rm{;}}\\
	\gamma,\ \hbox{ \rm{if }$X$ \rm{is of Type 3.2}\rm{.}}\\
	\end{array}
	\right.\nonumber
	\end{eqnarray}
Hence, $$N_5+2(N_2+N_4^{(1)}+N_6)+3(N_3^{(1)}+N_3^{(3)})+4N_4^{(2)}+\gamma N_3^{(2)}$$ $$\leq\sum_{X\in{\cal F}}|{\rm{supp}}(\Delta X)\cap A_{e\cdot}|\leq|A_{e\cdot}\setminus\{(0,0)\}|=|A_{e\cdot}|-1=\frac{mn}{2}-1.$$
Note that the second inequality comes from the cross-correlation parameter $\lambda_c=1$. \qed

\begin{Lemma}\label{1-13}
For any $m\equiv 0\ ({\rm mod }\ 4)$ and $\lambda_a\in\{2,3\}$,
\begin{eqnarray}\label{e-13}\nonumber
\Theta(m,n,3,\lambda_a,1)\leq \left\lfloor\frac{1}{8}\left(\frac{3mn}{2}+\eta+2\rho+4\omega-2+(3-\gamma)\varepsilon\right)\right\rfloor.
\end{eqnarray}
\end{Lemma}

\proof By (\ref{mn-1})+$4\times$(\ref{N2})+$2\times$(\ref{N3andN5})+(\ref{N41})+(\ref{mn/2-1}), we have
$$8N_2+12N_3^{(1)}+(5+\gamma)N_3^{(2)}+8N_3^{(3)}+8N_4^{(1)}+8N_4^{(2)}+8N_5+8N_6\leq \frac{3mn}{2}+\eta+2\rho+4\omega-2.$$
It follows that $$8(N_2+N_3^{(1)}+N_3^{(2)}+N_3^{(3)}+N_4^{(1)}+N_4^{(2)}+N_5+N_6)\leq \frac{3mn}{2}+\eta+2\rho+4\omega-2+(3-\gamma)N_3^{(2)}-4N_3^{(1)}.$$
By (\ref{N32}), $N_3^{(2)}\leq\varepsilon$. Due to $N_3^{(1)}\geq0$, we have
\begin{eqnarray}
8\Theta(m,n,3,\lambda_a,1)\leq\frac{3mn}{2}+\eta+2\rho+4\omega-2+(3-\gamma)\varepsilon.\nonumber
\end{eqnarray}
Then the conclusion follows immediately. \qed

\begin{Remark}\label{1-r-2}
Examining the proof of Lemma $\ref{1-13}$, we have that if
\begin{eqnarray}\nonumber \Theta(m,n,3,\lambda_a,1)=\frac{1}{8}\left(\frac{3mn}{2}+\eta+2\rho+4\omega-2+(3-\gamma)\varepsilon\right),
\end{eqnarray}
where $m\equiv 0\ ({\rm mod }\ 4)$ and $\lambda_a\in\{2,3\}$, then $N_3^{(1)}=0$, and the equalities must hold in \eqref{mn-1}, $\eqref{N2}$, $\eqref{N3andN5}$, $\eqref{N41}$ and $\eqref{mn/2-1}$. Especially, when $3-\gamma\neq 0$, i.e., $m\equiv 4\ ({\rm mod }\ 8)$, $N_3^{(2)}=\varepsilon$.
\end{Remark}

Input the exact values of $\omega$, $\varepsilon$, $\rho$, $\eta$, $\gamma$ to Lemma \ref{1-13}, and combine with Lemmas \ref{odd-m} and \ref{1-8}. We get an explicit upper bound of $\Theta(m,n,3,\lambda_a,1)$ for any positive integers $m$ and $n$.

\begin{Theorem}\label{1-main}
Let $\lambda_a\in\{2,3\}$. Then $\Theta(m,n,3,\lambda_a,1)\leq$
\begin{eqnarray}\label{e-10}\nonumber
\left\{
	\begin{array}{ll}
    \left\lfloor\frac{ mn+2\omega}{4}\right\rfloor,
	& \hbox{$mn\equiv 1,2,3\ ({\rm mod }\ 4)$\rm{;}}\\[0.4em]
	
	\left\lfloor\frac{7mn+16\omega}{32}\right\rfloor,&
	\hbox{$mn\equiv 0\ ({\rm mod }\ 8)$ \rm{and} $\gcd(m,n,2)=1$\rm{;}}\\[0.4em]
	
	\left\lfloor\frac{7mn+4+16\omega}{32}\right\rfloor,&
	\hbox{$mn\equiv 4\ ({\rm mod }\ 8)$ \rm{and} $\gcd(m,n,2)=1$\rm{;}}\\[0.4em]
	
	\left\lfloor\frac{5mn+4+8\omega}{24}\right\rfloor,
	& \hbox{$mn\equiv 4\ ({\rm mod }\ 8)$ \rm{and} $\gcd(m,n,2)=2$\rm{;}}\\[0.4em]

	\left\lfloor\frac{13mn+40+32\omega}{64}\right\rfloor,&
	\hbox{$mn\equiv 8\ ({\rm mod }\ 16),$  $\gcd(m,n,2)=2$, \rm{and} }\\
    &\hbox{$(m,n)\not\in\{(2,12),(4,6),(6,4),(6,12),(12,2),(12,6)\}$ \rm{when} $\lambda_a=3$\rm{;}}\\[0.4em]
	
	\left\lfloor\frac{5mn+8+8\omega}{24}\right\rfloor, &\hbox{$(m,n)\in\{(2,12),(4,6),(6,4),(6,12),(12,2),(12,6)\}$ \rm{and} $\lambda_a=3$\rm{;}} \\[0.4em]
	
	\left\lfloor\frac{13mn+32+32\omega}{64}\right\rfloor,
	&\hbox{$mn\equiv 0\ ({\rm mod }\ 32)$ \rm{and} $\gcd(m,n,2)=2$\rm{,} \rm{or }}\\
	&\hbox{$mn\equiv 16\ ({\rm mod }\ 32)$ \rm{and} $\gcd(m,n,4)=2$\rm{;}}\\[0.4em]
	
	\left\lfloor\frac{13mn+48+32\omega}{64}\right\rfloor,& \hbox{$mn\equiv 16\ ({\rm mod }\ 32)$ \rm{and}
	$\gcd(m,n,4)=4$\rm{.}}\\

	\end{array}
	\right.
	\end{eqnarray}
\end{Theorem}

\proof Note that $\Theta(m,n,3,\lambda_a,1)=\Theta(n,m,3,\lambda_a,1)$.

For $mn\equiv 1,2,3\ ({\rm mod }\ 4)$, at least one of $m$ and $n$ is odd. Then the conclusion follows from Lemma \ref{odd-m}.

For $mn\equiv 0\ ({\rm mod }\ 4)$ and $\gcd(m,n,2)=1$, w.l.o.g., assume that $n\equiv 1\ ({\rm mod }\ 2)$ and $m\equiv 0\ ({\rm mod }\ 4)$. Apply Lemma \ref{1-13} with $\varepsilon=\rho=1$, $(\eta,\gamma)=(\frac{mn}{4},3)$ and $(\frac{mn}{4}-1,1)$ when $m\equiv 0\ ({\rm mod }\ 8)$ and $m\equiv 4\ ({\rm mod }\ 8)$, respectively. Then we have
\begin{eqnarray}\label{e-10}\nonumber
	\Theta(m,n,3,\lambda_a,1)\leq\left\{
\begin{array}{ll}
	
	\left\lfloor\frac{7mn+16\omega}{32}\right\rfloor,&
	\hbox{\rm{if} $mn\equiv 0\ ({\rm mod }\ 8)$ \rm{and} $\gcd(m,n,2)=1$\rm{;}}\\[0.2em]

	\left\lfloor\frac{7mn+4+16\omega}{32}\right\rfloor,&
	\hbox{\rm{if} $mn\equiv 4\ ({\rm mod }\ 8)$ \rm{and} $\gcd(m,n,2)=1$\rm{.}}\\
\end{array}
\right.
\end{eqnarray}

For $mn\equiv 4\ ({\rm mod }\ 8)$ and $\gcd(m,n,2)=2$, we have $m,n\equiv 2\ ({\rm mod }\ 4)$. Then the conclusion follows from Lemma \ref{1-8}.

For $mn\equiv 8\ ({\rm mod }\ 16)$ and $\gcd(m,n,2)=2$, w.l.o.g., assume that $m\equiv 2\ ({\rm mod }\ 4)$ and $n\equiv 4\ ({\rm mod }\ 8)$.
By Lemma \ref{1-8}, we have
$$\Theta(m,n,3,\lambda_a,1)\leq\left\lfloor\frac{5mn+8+8\omega}{24}\right\rfloor:=U_1.$$
W.l.o.g., we can also assume that $m\equiv 4\ ({\rm mod }\ 8)$ and $n\equiv 2\ ({\rm mod }\ 4)$.
Applying Lemma \ref{1-13} with $\varepsilon=\gamma=1$, $\rho=3$ and $\eta=\frac{mn}{8}-1$, we have
$$\Theta(m,n,3,\lambda_a,1)\leq\left\lfloor\frac{13mn+40+32\omega}{64}\right\rfloor:=U_2.$$
It follows that
$\Theta(m,n,3,\lambda_a,1)\leq\min\{U_1,U_2\}.$ Comparing the values of $U_1$ and $U_2$, we have
\begin{eqnarray}\nonumber
\min\{U_1,U_2\}=\left\{
	\begin{array}{ll}
	U_2, &\hbox{\rm{if} $m\equiv 4\ ({\rm mod }\ 8),$  $n\equiv 2\ ({\rm mod }\ 4)${\rm ,}  \rm{and} } \\
	&\hbox{$(m,n)\not\in\{(4,6),(12,2),(12,6)\}$ \rm{when} $\lambda_a=3$\rm{;}}\\[0.2em]
	U_1, &\hbox{\rm{if} $(m,n)\in\{(4,6),(12,2),(12,6)\}$ \rm{and} $\lambda_a=3$\rm{.}} \\
	\end{array}
	\right.
\end{eqnarray}

For $mn\equiv 0\ ({\rm mod }\ 32)$ and $\gcd(m,n,2)=2$, w.l.o.g., assume that $m\equiv 0\ ({\rm mod }\ 8)$ and $n\equiv 0\ ({\rm mod }\ 2)$.
By Lemma \ref{1-13} with $\rho=\gamma=3$ and $\eta=\frac{mn}{8}$, we have
$$\Theta(m,n,3,\lambda_a,1)\leq \left\lfloor\frac{13mn+32+32\omega}{64}\right\rfloor:=U_3.$$
W.l.o.g., we can also assume that $m\equiv 0\ ({\rm mod }\ 2)$ and $n\equiv 0\ ({\rm mod }\ 8)$. Applying Lemma \ref{1-8} with $m\equiv 2\ ({\rm mod }\ 4)$, and Lemma \ref{1-13} with $\varepsilon=2$, $\rho=3$, $(\eta,\gamma)=(\frac{mn}{8}-2,1)$ and $(\frac{mn}{8},3)$ for $m\equiv 4\ ({\rm mod }\ 8)$ and $m\equiv 0\ ({\rm mod }\ 8)$, respectively, we have
\begin{eqnarray}\nonumber
\Theta(m,n,3,\lambda_a,1)\leq\left\{
	\begin{array}{ll}
	U_1, &\hbox{  \rm{if} $m\equiv 2\ ({\rm mod }\ 4)${\rm ;}} \\[0.2em]
	\left\lfloor\frac{13mn+48+32\omega}{64}\right\rfloor:=U_4, &\hbox{  \rm{if} $m\equiv 4\ ({\rm mod }\ 8)${\rm ;}} \\[0.2em]
	U_3, &\hbox{  \rm{if} $m\equiv 0\ ({\rm mod }\ 8)$.}
	\end{array}
	\right.
\end{eqnarray}
Therefore,
\begin{eqnarray}\nonumber
\Theta(m,n,3,\lambda_a,1)\leq\left\{
	\begin{array}{ll}
	\min\{U_3,U_1\}=U_3, &\hbox{  \rm{if} $m\equiv 2\ ({\rm mod }\ 4)${\rm ;}} \\
	\min\{U_3,U_4\}=U_3, &\hbox{  \rm{if} $m\equiv 4\ ({\rm mod }\ 8)$\rm{;}} \\
	U_3, &\hbox{  \rm{if} $m\equiv 0\ ({\rm mod }\ 8)$\rm{.}} \\
	\end{array}
	\right.
\end{eqnarray}

For $mn\equiv 16\ ({\rm mod }\ 32)$ and $\gcd(m,n,2)=2$, we consider two subcases.
If $mn\equiv 16\ ({\rm mod }\ 32)$ and $\gcd(m,n,4)=2$, then assume that $m\equiv 2\ ({\rm mod }\ 4)$ and $n\equiv 8\ ({\rm mod }\ 16)$. By Lemma \ref{1-8}, we have
$\Theta(m,n,3,\lambda_a,1)\leq U_1.$ We can also assume that $m\equiv 8\ ({\rm mod }\ 16)$ and $n\equiv 2\ ({\rm mod }\ 4)$. Then applying Lemma \ref{1-13} with $\rho=\gamma=3$ and $\eta=\frac{mn}{8}$, we have $\Theta(m,n,3,\lambda_a,1)\leq U_3.$ Therefore, we get $\Theta(m,n,3,\lambda_a,1)\leq\min\{U_1,U_3\}=U_3.$ If $mn\equiv 16\ ({\rm mod }\ 32)$ and $\gcd(m,n,4)=4$, which implies $m,n\equiv 4\ ({\rm mod }\ 8)$, then by Lemma \ref{1-13} with $\varepsilon=2$, $\rho=3$, $\eta=\frac{mn}{8}-2$ and $\gamma=1$, we have $\Theta(m,n,3,\lambda_a,1)\leq U_4$. \qed

\subsection{\label{321} Improved upper bound for two subclasses when $\lambda_a=2$}

A codeword of Type $4.2$ is of the form $\{(0,0),(a,b),(2a,2b)\}$, where $a\equiv 0\ ({\rm mod }\ 2)$. All codewords of Type $4.2$ can be divided into the following three types:
\begin{itemize}
\item[]\textbf{Type 4.2.1}: $a\equiv 2\ ({\rm mod }\ 4)$ and $b\equiv 0\ ({\rm mod }\ 2)$;
\item[]\textbf{Type 4.2.2}: $a\equiv 0\ ({\rm mod }\ 4)$ and $b\equiv 0\ ({\rm mod }\ 2)$;
\item[]\textbf{Type 4.2.3}: $a\equiv 0\ ({\rm mod }\ 2)$ and $b\equiv 1\ ({\rm mod }\ 2)$.
\end{itemize}
\noindent Let $N_4^{(2,1)}$, $N_4^{(2,2)}$ and $N_4^{(2,3)}$ denote the number of codewords in $\cal{F}$ of Types $4.2.1$, $4.2.2$ and $4.2.3$, respectively.
Then, $$N_4^{(2)}=N_4^{(2,1)}+N_4^{(2,2)}+N_4^{(2,3)}.$$

% The following lemma is useful to improve the upper bound of $\Theta(m,n,3,2,1)$.

\begin{Lemma}\label{2-0} If $\lambda_a=2$ and $2N_4^{(1)}=\eta$,
	then
	\begin{eqnarray}\nonumber
	4N_4^{(2)}\leq\left\{
	\begin{array}{ll}
\frac{mn}{4}-2\xi-1,&
\ \hbox{ \rm{if }$m\equiv 4\ ({\rm mod }\ 8)$ \rm{and} $n\equiv 1\ ({\rm mod }\ 2)$\rm{;}}\\[0.2em]
\frac{3mn}{16}-3,&\ \hbox{ \rm{if }$m\equiv 4,20\ ({\rm mod }\ 24)$ \rm{and }$n\equiv 4\ ({\rm mod }\ 8)$\rm{;}}\\[0.2em]
\frac{3mn}{16}-5,&\ \hbox{ \rm{if }$m\equiv 12\ ({\rm mod }\ 24)$ \rm{and }$n\equiv 4,20\ ({\rm mod }\ 24)$\rm{;}}\\[0.2em]
\frac{3mn}{16}-7,&\ \hbox{ \rm{if }$m,n\equiv 12\ ({\rm mod }\ 24)$\rm{;}}\\[0.2em]
\frac{3mn}{16}-6,&
\ \hbox{ \rm{if }$m\equiv 0\ ({\rm mod }\ 8)$ \rm{and} $n\equiv 4\ ({\rm mod }\ 8)$\rm{.}}\\
\end{array}
	\right.
	\end{eqnarray}
Furthermore, when $m\equiv 0\ ({\rm mod }\ 8)$ and $n\equiv 4\ ({\rm mod }\ 8)$, if $4N_4^{(2)}=\frac{3mn}{16}-6$, then $2N_4^{(2,3)}=\frac{mn}{16}-2$.
\end{Lemma}

\proof \noindent\textbf{Case 1}: $m\equiv 4\ ({\rm mod }\ 8)$ and $n\equiv 1\ ({\rm mod }\ 2)$.

Due to $A_{s\cdot}\setminus\Omega_A(2)=A_{s\cdot}\setminus\{(\frac{m}{2},0)\}$, we have $|A_{s\cdot}\setminus\Omega_A(2)|=\frac{m}{4}\times n-1=\eta$.
Since each codeword of Type $4.1$ contributes exactly two differences in $A_{s\cdot}\setminus\Omega_A(2)$, the condition $2N_4^{(1)}=\eta$ implies that every element in $A_{s\cdot}\setminus\Omega_A(2)$ is used as a difference of some codeword of Type $4.1$. Hence for each codeword $\{(0,0),(a,b),(2a,2b)\}$ of Type $4.2$, we have $a\equiv 0\ ({\rm mod }\ 4)$. It follows that each codeword of Type $4.2$ contributes four differences in $A_{d\cdot}$. Since $m\equiv 4\ ({\rm mod }\ 8)$ and $n\equiv 1\ ({\rm mod }\ 2)$, it is readily checked that for any codeword $X$ of Type $4.2$, we have ${\rm supp}(\Delta X)\cap\Omega_A(3)=\emptyset$ (otherwise, either $\lambda_a=3$ or $4\mid\frac{m}{6}$). Due to
$$\begin{array}{lll}
|A_{d\cdot}\cap\Omega_A(3)|=\left\{
\begin{array}{ll}
1,\ \hbox{ \rm{if }$3\nmid mn$\rm{;}}\\
3,\ \hbox{ \rm{if }$3\mid mn$ {\rm and} $\gcd(m,n,3)=1$\rm{;}}\\
9,\ \hbox{ \rm{if} $\gcd(m,n,3)=3$,}\\
\end{array}
\right.
\end{array}$$
we have $|A_{d\cdot}\cap\Omega_A(3)|=2\xi+1$. Therefore,
$$4N_4^{(2)}\leq |A_{d\cdot}\setminus\Omega_A(3)|=|A_{d\cdot}|-|A_{d\cdot}\cap\Omega_A(3)|
=\frac{m}{4}\times n-(2\xi+1)=\frac{mn}{4}-2\xi-1.
$$

\noindent\textbf{Case 2}: $m\equiv 0\ ({\rm mod }\ 4)$ and $n\equiv 4\ ({\rm mod }\ 8)$.

For $m\equiv 0\ ({\rm mod }\ 8)$ and $n\equiv 4\ ({\rm mod }\ 8)$, $A_{se}\setminus\Omega_A(2)=A_{se}$. Hence $|A_{se}\setminus\Omega_A(2)|=\frac{m}{4}\times\frac{n}{2}=\frac{mn}{8}=\eta$. For $m\equiv 4\ ({\rm mod }\ 8)$ and $n\equiv 4\ ({\rm mod }\ 8)$, $A_{se}\setminus\Omega_A(2)=A_{se}\setminus\{(\frac{m}{2},0),(\frac{m}{2},\frac{n}{2})\}$. Hence $|A_{se}\setminus\Omega_A(2)|=\frac{m}{4}\times\frac{n}{2}-2=\frac{mn}{8}-2=\eta$.

Let $\{(0,0),(a,b),(2a,2b)\}$ be a codeword of Type $4.1$, where $a\equiv 1\ ({\rm mod }\ 2)$. It contributes two differences $(2a,2b)$ and $(-2a,-2b)$ in $A_{se}\setminus\Omega_A(2)$. Due to $|A_{se}\setminus\Omega_A(2)|=\eta$, the condition $2N_4^{(1)}=\eta$ implies that every element in $A_{se}\setminus\Omega_A(2)$ is used as a difference of some codeword of Type $4.1$. It follows that
\begin{eqnarray}\label{e0}
N_4^{(2,1)}=0.
\end{eqnarray}

Due to
$$\begin{array}{lll}
A_{ds}\cap\Omega_A(2)=\left\{
\begin{array}{ll}
\{(0,\frac{n}{2})\},&
\ \hbox{ \rm{if }$m\equiv 4\ ({\rm mod }\ 8)$\rm{;}}\\[0.2em]
\{(0,\frac{n}{2}),(\frac{m}{2},\frac{n}{2})\},&
\ \hbox{ \rm{if }$m\equiv 0\ ({\rm mod }\ 8)$\rm{,}}\\
\end{array}
\right.
\end{array}$$
we have
$$
|A_{ds}\cap\Omega_A(2)|=\left\{
\begin{array}{ll}
1,
\ \hbox{ \rm{if }$m\equiv 4\ ({\rm mod }\ 8)$\rm{;}}\\
2,
\ \hbox{ \rm{if }$m\equiv 0\ ({\rm mod }\ 8)$\rm{.}}\\
\end{array}
\right.
$$
For any codeword $X$ of Type $4.2.3$, $|{\rm supp}(\Delta X)\cap A_{ds}|=|{\rm supp}(\Delta X)\cap A_{de}|=2$. By the definition of ${\cal T}_1$, each codeword $X$ of Type $4.2$ satisfies ${\rm supp}(\Delta X)\cap\Omega_A(2)=\emptyset$. Therefore,
\begin{eqnarray}\label{e1}
2N_4^{(2,3)}\leq |A_{ds}\setminus\Omega_A(2)|=|A_{ds}|-|A_{ds}\cap\Omega_A(2)|=\left\{
\begin{array}{ll}
\frac{mn}{16}-1,
\ \hbox{ \rm{if }$m\equiv 4\ ({\rm mod }\ 8)$\rm{;}}\\[0.2em]
\frac{mn}{16}-2,
\ \hbox{ \rm{if }$m\equiv 0\ ({\rm mod }\ 8)$\rm{.}}\\
\end{array}
\right.
\end{eqnarray}

Due to
$$\begin{array}{lll}
A_{de}\cap\Omega_A(2)=\left\{
\begin{array}{ll}
\{(0,0),(0,\frac{n}{2})\},&
\ \hbox{ \rm{if }$m\equiv 4\ ({\rm mod }\ 8)$\rm{;}}\\[0.2em]
\{(0,0),(0,\frac{n}{2}),(\frac{m}{2},0),(\frac{m}{2},\frac{n}{2})\},&
\ \hbox{ \rm{if }$m\equiv 0\ ({\rm mod }\ 8)$,}\\
\end{array}
\right.\\
\end{array}$$
we have
$$
|A_{de}\cap\Omega_A(2)|=\left\{
\begin{array}{ll}
2,
\ \hbox{ \rm{if }$m\equiv 4\ ({\rm mod }\ 8)$\rm{;}}\\
4,
\ \hbox{ \rm{if }$m\equiv 0\ ({\rm mod }\ 8)$\rm{.}}\\
\end{array}
\right.
$$
For any codeword $X$ of Type $4.2.2$, we have $|{\rm supp}(\Delta X)\cap A_{de}|=4$. Therefore, when
$m\equiv 0\ ({\rm mod }\ 8)$,
\begin{eqnarray}\label{e2}
4N_4^{(2,2)}+2N_4^{(2,3)}\leq|A_{de}\setminus\Omega_A(2)|=|A_{de}|-|A_{de}\cap\Omega_A(2)|=\frac{mn}{8}-4.
\end{eqnarray}
Let $W=\left\{\pm(0,\frac{n}{3}),\pm(\frac{2m}{3},\frac{n}{3})\right\}$. Then $W\subset \Omega_A(3)\subset A_{de}$. Hence,
$$\left| A_{de}\cap (\Omega_A(3)\setminus W)\right|=\left|\Omega_A(3)\setminus W\right|=\left\{
\begin{array}{lll}
1, & {\rm if }\ m\equiv 4,20\ ({\rm mod }\ 24)\ {\rm and}\ n\equiv 4\ ({\rm mod }\ 8);\\[0.2em]
3, & {\rm if }\ m\equiv 12\ ({\rm mod }\ 24)\ {\rm and}\ n\equiv 4,20\ ({\rm mod }\ 24);\\[0.2em]
5, & {\rm if }\ m,n\equiv 12\ ({\rm mod }\ 24).\\
\end{array}
\right.
$$
When $m,n\equiv 4\ ({\rm mod }\ 8)$, it is readily checked that for any codeword $X$ of Type $4.2.2$ or Type $4.2.3$, we have ${\rm supp}(\Delta X)\cap\Omega_A(3)\subseteq W$ (otherwise, either $\lambda_a=3$ or $4\mid\frac{m}{6}$). Therefore,
\begin{eqnarray}\label{e3}
4N_4^{(2,2)}+2N_4^{(2,3)}&\leq&|A_{de}\setminus(\Omega_A(2)\cup(\Omega_A(3)\setminus W))| \nonumber \\\nonumber
&=&|A_{de}|-|A_{de}\cap\Omega_A(2)|-|\Omega_A(3)\setminus W|+1\\
&=&\left\{
\begin{array}{lll}
\frac{mn}{8}-2,&\ \hbox{ \rm{if }$m\equiv 4,20\ ({\rm mod }\ 24)$ \rm{and }$n\equiv 4\ ({\rm mod }\ 8)$\rm{;}}\\[0.2em]
\frac{mn}{8}-4,&\ \hbox{ \rm{if }$m\equiv 12\ ({\rm mod }\ 24)$ \rm{and }$n\equiv 4,20\ ({\rm mod }\ 24)$\rm{;}}\\[0.2em]
\frac{mn}{8}-6,&\ \hbox{ \rm{if }$m,n\equiv 12\ ({\rm mod }\ 24)$\rm{.}}\\
\end{array}
\right.
\end{eqnarray}

\noindent By $4\times\eqref{e0}+\eqref{e1}+\eqref{e2}$ and $4\times\eqref{e0}+\eqref{e1}+\eqref{e3}$, we obtain
$$4N_4^{(2)}=4N_4^{(2,2)}+4N_4^{(2,3)}\leq\left\{
\begin{array}{lll}
\frac{3mn}{16}-3,&\ \hbox{ \rm{if }$m\equiv 4,20\ ({\rm mod }\ 24)$ \rm{and }$n\equiv 4\ ({\rm mod }\ 8)$\rm{;}}\\[0.2em]
\frac{3mn}{16}-5,&\ \hbox{ \rm{if }$m\equiv 12\ ({\rm mod }\ 24)$ \rm{and }$n\equiv 4,20\ ({\rm mod }\ 24)$\rm{;}}\\[0.2em]
\frac{3mn}{16}-7,&\ \hbox{ \rm{if }$m,n\equiv 12\ ({\rm mod }\ 24)$\rm{;}}\\[0.2em]
\frac{3mn}{16}-6,&\ \hbox{ \rm{if }$m\equiv 0\ ({\rm mod }\ 8)$ \rm{and }$n\equiv 4\ ({\rm mod }\ 8)$\rm{.}}\\
\end{array}
\right.
$$
Furthermore, when $m\equiv 0\ ({\rm mod }\ 8)$ and $n\equiv 4\ ({\rm mod }\ 8)$, if $4N_4^{(2)}=\frac{3mn}{16}-6$, then the equality holds in \eqref{e1}, i.e., $2N_4^{(2,3)}=\frac{mn}{16}-2$. \qed

\begin{Lemma}\label{2-eg-3} Let $m,n\equiv 4\ ({\rm mod }\ 16)$ or $m,n\equiv 12\ ({\rm mod }\ 16)$ be positive integers such that $3\mid m$. Then $$\Theta(m,n,3,2,1)\leq \frac{13mn-16}{64}.$$
\end{Lemma}

\proof For $m,n\equiv 4\ ({\rm mod }\ 16)$, or $m,n\equiv 12\ ({\rm mod }\ 16)$, by Theorem \ref{1-main} with $\omega=0$, $\Theta(m,n,3,2,1)\leq \frac{13mn+48}{64}$. Assume that $\Theta(m,n,3,2,1)=\frac{13mn+48}{64}$. By Remark \ref{1-r-2}, $N_3^{(1)}=0$, $N_3^{(2)}=\varepsilon=2$ and the equalities hold in (\ref{mn-1}), (\ref{N3andN5}), (\ref{N41}) and (\ref{mn/2-1}). Then by (\ref{N3andN5}), $N_3^{(3)}+N_5=1$. By (\ref{N41}), $N_4^{(1)}=\frac{mn}{16}-1$.  By $\eqref{mn-1}-\eqref{mn/2-1}$, we have
$$2N_3^{(2)}+2N_4^{(1)}+4N_5+4N_6=\frac{mn}{2},$$
which yields
$$N_6= \frac{3mn-16}{32}-N_5\leq\frac{3mn-16}{32}.$$
Note that $N_2=0$. Therefore,
$$4N_4^{(2)}=4\left(\Theta(m,n,3,2,1)-N_4^{(1)}-N_3^{(2)}-(N_3^{(3)}+N_5)-N_6\right)=\frac{9mn}{16}-5-4N_6\geq\frac{3mn}{16}-3.$$
However, the condition $3\mid m$ implies that $m\equiv 12\ ({\rm mod }\ 24)$. Hence by Lemma \ref{2-0},
$4N_4^{(2)}\leq\frac{3mn}{16}-7$ or $\frac{3mn}{16}-5$ according to whether $n$ is divided by $3$ or not,
a contradiction. \qed

A codeword of Type $4.1$ is of the form $\{(0,0),(a,b),(2a,2b)\}$, where $a\equiv 1\ ({\rm mod }\ 2)$. All the codewords of Type $4.1$ can be divided into the following two types according to the parity of $b$:
\begin{itemize}
\item[]\textbf{Type 4.1.1}: $a,b\equiv 1\ ({\rm mod }\ 2)$;
\item[]\textbf{Type 4.1.2}: $a\equiv 1\ ({\rm mod }\ 2)$ and $b\equiv 0\ ({\rm mod }\ 2)$;
\end{itemize}
\noindent Let $N_4^{(1,1)}$ and $N_4^{(1,2)}$ denote the number of codewords in $\cal{F}$ of Types $4.1.1$ and $4.1.2$, respectively.
Then, $$N_4^{(1)}=N_4^{(1,1)}+N_4^{(1,2)}.$$

\begin{Lemma}\label{N4112}
Let $m\equiv 0\ ({\rm mod }\ 8)$ and $n\equiv 0\ ({\rm mod }\ 4)$. Assume that $2N_4^{(1)}=\frac{mn}{8}$.
Then $N_4^{(1,1)}$ $=N_4^{(1,2)}=\frac{mn}{32}.$ Furthermore, if $m\equiv 8\ ({\rm mod }\ 16)$, then $N_3^{(2)}\leq1$, and for any codeword $X$ of Type $3.2$, ${\rm supp}(\Delta X)=\left\{\pm\left(\frac{m}{4},\frac{n}{4}\right),\left(\frac{m}{2},\frac{n}{2}\right)\right\}$
or $\left\{\pm\left(\frac{m}{4},\frac{3n}{4}\right),\left(\frac{m}{2},\frac{n}{2}\right)\right\}$.
\end{Lemma}

\proof Let $\{(0,0),(a,b),(2a,2b)\}$ with $a\equiv 1\ ({\rm mod }\ 2)$ be a codeword of Type $4.1$. It contributes two differences $(2a,2b)$ and $(-2a,-2b)$ in $A_{se}$. Due to $|A_{se}|=\frac{m}{4}\times\frac{n}{2}=\frac{mn}{8}$, the condition $2N_4^{(1)}=\frac{mn}{8}$ implies that every element in $A_{se}$ is used as a difference of some codeword of Type $4.1$. Therefore, $n\equiv 0\ ({\rm mod }\ 4)$ yields $N_4^{(1,1)}=N_4^{(1,2)}=\frac{1}{2}N_4^{(1)}=\frac{mn}{32}$. Furthermore, when $m\equiv 8\ ({\rm mod }\ 16)$,
$\left\{\pm\left(\frac{m}{4},0\right), \pm\left(\frac{m}{4},\frac{n}{2}\right)\right\}\subset A_{se}$. So
$\pm(\frac{m}{4},0)$ and $\pm(\frac{m}{4},\frac{n}{2})$ are used as differences of some codewords of Type $4.1$, and cannot be produced by other types of codewords. Therefore, for any codeword $X$ of Type $3.2$, by Remark \ref{for type3.1.2.3}(2),
${\rm supp}(\Delta X)=\left\{\pm\left(\frac{m}{4},\frac{n}{4}\right),\left(\frac{m}{2},\frac{n}{2}\right)\right\}$
or $\left\{\pm\left(\frac{m}{4},\frac{3n}{4}\right),\left(\frac{m}{2},\frac{n}{2}\right)\right\}$.
The two sets share a common element $\left(\frac{m}{2},\frac{n}{2}\right)$, so $N_3^{(2)}\leq1$. \qed

\begin{Lemma}\label{2-eg-2} For any $m\equiv 8\ ({\rm mod }\ 16)$ and $n\equiv 4\ ({\rm mod }\ 8)$, $$\Theta(m,n,3,2,1)\leq \frac{13mn-32}{64}.$$
\end{Lemma}

\proof For $m\equiv 8\ ({\rm mod }\ 16)$ and $n\equiv 4\ ({\rm mod }\ 8)$, applying Theorem \ref{1-main} with $\omega=0$, we have $\Theta(m,n,3,2,1)\leq \frac{13mn+32}{64}$. Assume that $\Theta(m,n,3,2,1)=\frac{13mn+32}{64}$. By Remark \ref{1-r-2}, $N_3^{(1)}=0$ and the equalities hold in \eqref{mn-1}, (\ref{N3andN5}), (\ref{N41}) and (\ref{mn/2-1}). By (\ref{N41}), $2N_4^{(1)}=\frac{mn}{8}$. It follows that $N_3^{(2)}\leq1$ by Lemma \ref{N4112}. Note that $N_3^{(3)}\leq1$ by (\ref{N33}). Thus, by (\ref{N3andN5}), $N_3^{(2)}+N_3^{(3)}+N_5=3$ yields $N_5\geq1$. By (\ref{mn-1})$-$(\ref{mn/2-1}), we have $2N_4^{(1)}+4N_5+4N_6=\frac{mn}{2}$, which yields
$$N_6= \frac{3mn}{32}-N_5\leq\frac{3mn}{32}-1.$$
Note that $N_2=0$ by Lemma \ref{1-1}. Therefore,
$$N_4^{(2)}=\Theta(m,n,3,2,1)-N_4^{(1)}-(N_3^{(2)}+N_3^{(3)}+N_5)-N_6=\frac{9mn-160}{64}-N_6\geq\frac{3mn-96}{64}.$$
By Lemma \ref{2-0}, $N_4^{(2)}\leq\frac{3mn-96}{64}.$ Thus $N_4^{(2)}=\frac{3mn-96}{64}.$ It follows that
$N_3^{(2)}=N_3^{(3)}=N_5=1$ and $N_6=\frac{3mn}{32}-1.$

%Recall that the equality holds in (\ref{mn-1}), that is, the difference leave of ${\cal F}$ is empty.
Recall that the equality holds in (\ref{mn-1}). It follows that the elements in $A_{\cdot o}$ are used up as the differences of codewords of Types $3.2$, $3.3$, $4.1.1$, $4.2.3$, $5$ and $6$. Since $N_3^{(2)}=1$, by Lemma \ref{N4112}, the unique codeword $X$ of Type $3.2$ satisfies $|{\rm supp}(\Delta X)\cap A_{\cdot o}|=2$. Since $N_3^{(3)}=1$, the unique codeword $X$ of Type $3.3$ also satisfies $|{\rm supp}(\Delta X)\cap A_{\cdot o}|=2$ by Remark \ref{for type3.1.2.3}(3). Due to $2N_4^{(1)}=\frac{mn}{8}$, by Lemma \ref{N4112}, there are $N_4^{(1,1)}=\frac{mn}{32}$ codewords of Type $4.1.1$, and by Lemma \ref{2-0}, there are $N_4^{(2,3)}=\frac{mn}{32}-1$ codewords of Type $4.2.3$. Each codeword of Type $4.1.1$ (resp. of Type $4.2.3$) contributes two distinct differences in $A_{\cdot o}$.
Let $Q$ denote the number of distinct differences in $A_{\cdot o}$ from all codewords of Type $5$ and of Type $6$. Then it is readily checked that $Q\equiv 0\ ({\rm mod }\ 4)$. However, since $|A_{\cdot o}|=m\times\frac{n}{2}=\frac{mn}{2}$, we get $$2+2+\frac{mn}{32}\times2+(\frac{mn}{32}-1)\times2+Q=\frac{mn}{2},$$
which yields $Q=\frac{3mn}{8}-2\equiv 2\ ({\rm mod }\ 4)$, a contradiction. \qed

\subsection{\label{321} Sporadic values}

\begin{Lemma}\label{1-eg-2}
$\Theta(4,2,3,\lambda_a,1)\leq 1$ for any $\lambda_a\in\{2,3\}$.
\end{Lemma}

\proof By (\ref{N31}), (\ref{N32}), (\ref{N33}) and (\ref{N2andN4}), $N_2=N_3^{(3)}=N_4=0$, $N_3^{(1)}\leq1$ and $N_3^{(2)}\leq1$.
By (\ref{mn-1}), we have
\begin{eqnarray}\label{4and2}
3N_3^{(1)}+3N_3^{(2)}+5N_5+6N_6\leq7.
\end{eqnarray}
By Theorem \ref{1-main}, $\Theta(4,2,3,\lambda_a,1)\leq 2$. Assume that $\Theta(4,2,3,\lambda_a,1)$ $=2$, that is, $N_3^{(1)}+N_3^{(2)}+N_5+N_6=2$. It follows from (\ref{4and2}) that $N_3^{(1)}=N_3^{(2)}=1$ and $N_5=N_6=0$. Thus $3N_3^{(1)}+N_3^{(2)}+N_3^{(3)}+N_5=4$. It contradicts with (\ref{N3andN5}). \qed

\begin{Lemma}\label{1-eg-1}
$\Theta(12,3,3,3,1)\leq 9.$
\end{Lemma}

\proof By Theorem \ref{1-main}, $\Theta(12,3,3,3,1)\leq10$. Assume that $\Theta(12,3,3,3,1)=10$. By Remark \ref{1-r-2}, $N_3^{(1)}=0$, $N_3^{(2)}=1$ and the equalities hold in (\ref{mn-1}), (\ref{N2}), (\ref{N3andN5}) and (\ref{N41}). By (\ref{N3andN5}), $N_3^{(3)}=N_5=0$. By (\ref{N2}) and (\ref{N41}), $N_2=4$ and $N_4^{(1)}=4$. By (\ref{mn-1}), $2N_2+3N_3+4N_4+5N_5+6N_6=35$, which yields $2N_4^{(2)}+3N_6=4$. It follows that $N_4^{(2)}=2$ and $N_6=0$. Therefore, $N_2+N_3+N_4+N_5+N_6=11$. It contradicts with $\Theta(12,3,3,3,1)=10$. \qed

\begin{Lemma}\label{2-eg-1}
$\Theta(12,3,3,2,1)\leq 7.$
\end{Lemma}
\proof The proof is subtly different from that of Lemma \ref{1-eg-1}. By Theorem \ref{1-main}, $\Theta(12,3,3,2,1)$ $\leq 8$. Assume that $\Theta(12,3,3,2,1)=8$. By Remark \ref{1-r-2}, $N_3^{(1)}=0$, $N_3^{(2)}=1$ and the equalities hold in (\ref{mn-1}), (\ref{N3andN5}), (\ref{N41}) and (\ref{mn/2-1}). By (\ref{N3andN5}), $N_3^{(3)}=N_5=0$. By (\ref{N2}) and (\ref{N41}), $N_2=0$ and $N_4^{(1)}=4$. Therefore, by (\ref{mn-1}), $2N_4^{(2)}+3N_6=8$, and by (\ref{mn/2-1}), $2N_4^{(2)}+N_6=4$. It follows that $N_6=2$ and $N_4^{(2)}=1$. However, by Lemma \ref{2-0}, $N_4^{(2)}=0$, a contradiction. \qed

\subsection{Proof of Theorem \ref{c-main-1}}

For $(m,n)\in\{(2,4),(4,2),(3,12),(12,3)\}$, the conclusion follows from Lemmas \ref{1-eg-2}, \ref{1-eg-1} and \ref{2-eg-1}. For $mn\equiv 32\ ({\rm mod }\ 64)$, $\gcd(m,n,8)=4$ and $\lambda_a=2$, exactly one of $m$ and $n$ is divided by $4$ but not by $8$. W.l.o.g., assume that $n\equiv 4\ ({\rm mod }\ 8)$ and $m\equiv 8\ ({\rm mod }\ 16)$. Then the conclusion follows from Lemma \ref{2-eg-2}. For $mn\equiv 144\ ({\rm mod }\ 192)$, $\gcd(m,n,4)=4$ and $\lambda_a=2$, we have $3\mid mn$ and $mn\equiv 16\ ({\rm mod }\ 64)$. W.l.o.g., assume that $3\mid m$, and $m,n\equiv 4\ ({\rm mod }\ 16)$ or $m,n\equiv 12\ ({\rm mod }\ 16)$. Then the conclusion follows by Lemma \ref{2-eg-3}. All the other cases follow from Theorem \ref{1-main}. \qed

\section{\label{add} Recursive constructions}

Let $\cal{F}$ be an $(m, n, k, \lambda_a, 1)$-OOSPC. Define the \emph{$($difference$)$ leave} of ${\cal F}$, briefly ${\rm DL}({\cal F})$, as the set of all nonzero elements in $\mathbb{Z}_m\times \mathbb{Z}_n$ which are not covered by $\Delta{\cal F}=\bigcup_{X\in {\cal F}}\textrm{supp}(\Delta X)$. ${\cal F}$ is said to be $(s, t)$-\emph{regular} if ${\rm DL}({\cal F})\cup \{(0,0)\}$ forms an additive subgroup $S\times T$ of $\mathbb{Z}_m\times \mathbb{Z}_n$, where $S$ and $T$ are, respectively, the additive subgroups of order $s$ in $\mathbb{Z}_m$ and order $t$ in $\mathbb{Z}_n$.
%The following construction is straightforward but very useful.
% A $(1\times 1)$-\emph{regular} $(m, n, k, \lambda_a, 1)$-OOSPC, whose leave is empty, is said to be {\em perfect}.

% The following construction is simple but very useful.

\begin{Construction}\label{Hole}{\rm{(Filling Construction)}}
Suppose that there exist
\begin{enumerate}
\item[{\rm (1)}]an $(s, t)$-regular $(m,n,k,\lambda_a,1)$-{\rm OOSPC} ${\cal{F}}_1$
with $b_1$ codewords;
\item[{\rm (2)}]an $(s,t,k,\lambda_a,1)$-{\rm OOSPC} ${\cal{F}}_2$ with $b_2$ codewords.
\end{enumerate}
Then there exists an $(m,n,k,\lambda_a,1)$-{\rm OOSPC} with $b_1+b_2$ codewords{.}
Furthermore, if the given $(s,t,k,\lambda_a,1)$-{\rm OOSPC} is $(g, h)$-regular, then the resulting $(m,n,k,\lambda_a,1)$-{\rm OOSPC} is $(g, h)$-regular.
\end{Construction}
\proof Let us interpret all codewords of ${\cal{F}}_2$ as codewords in $(\frac{m}{s}\mathbb Z_{m})\times (\frac{n}{t}\mathbb Z_n)$ and add them to the codewords of ${\cal{F}}_1$. We then get the desired $(m,n,k,\lambda_a,1)$-{\rm OOSPC} with $b_1+b_2$ codewords, whose leave is exactly DL$({\cal{F}}_2)$.
\qed

Let $G$ be an abelian group of order $v$.
A $(G,k,\lambda)$ \emph{difference matrix} (briefly, $(G,k,\lambda)$-DM) is a
$k\times \lambda v$ matrix $D=(d_{ij})$ with entries from $G$ such that for any
distinct rows $x$ and $y$, the multiset $\{d_{xi}-d_{yi}:\, 1\leq i\leq \lambda v\}$
contains each element of $G$ exactly $\lambda$ times.
If $G=\mathbb Z_v$, the difference matrix is called {\em cyclic} and denoted by a $(v,k,\lambda)$-CDM.

When $\lambda_{a}=1$, the notation $(s, t)$-regular $(m, n, k, 1, 1)$-OOSPC is simply written as $(s, t)$-regular $(m, n, k, 1)$-OOSPC.

\begin{Construction}{\rm \cite[Construction 3.3]{PC}} {\rm{(Inflation Construction)}}
\label{R-DM} Let $m,n$ and $v$ be positive integers. Suppose
that there exist
\begin{enumerate}
\item[{\rm (1)}] an $(s, t)$-regular $(m,n,k,1)$-{\rm OOSPC};
\item[{\rm (2)}] a $(v,k,1)$-{\rm CDM}.
\end{enumerate}
Then there exist an $(sv, t)$-regular $(mv,n,k,1)$-{\rm OOSPC} and an $(s, tv)$-regular $(m,nv,k,1)$-{\rm OOSPC}.
\end{Construction}

% In the rest of this paper, Lemma \ref{R-DM} is used many times and the existence of the needed $(v,k,1)$-CDM is guaranteed by the following Lemma \ref{CDM}. Hence, we will not repeatedly say that ``the needed $(v,k,1)$-CDM exists by Lemma \ref{CDM}".

\begin{Lemma}\label{CDM}\emph{\cite{CJC}}
Let $v$ and $k$ be positive integers such that $\gcd(v,(k-1)!)=1.$ Then there exists a $(v,k,1)$-{\rm CDM}.
\end{Lemma}

Let $m,n\equiv 2\ ({\rm mod}\ 4)$, and $H$ be the normal subgroup $\{(0,0),(0,\frac{n}{2}),(\frac{m}{2},0),(\frac{m}{2},\frac{n}{2})\}$ of $\mathbb{Z}_m\times \mathbb{Z}_n$. Then the quotient group $(\mathbb{Z}_m\times \mathbb{Z}_n)/H$ is isomorphic to $\mathbb{Z}_\frac{m}{2}\times \mathbb{Z}_\frac{n}{2}$. For $(x,y)\in \mathbb{Z}_\frac{m}{2}\times \mathbb{Z}_\frac{n}{2}$, let $D(x,y)=(x,y)+H$ be a coset of $H$ in $\mathbb{Z}_m\times \mathbb{Z}_n$, namely,
$$D(x,y)=\left\{(x,y),(x,y+\frac{n}{2}),(x+\frac{m}{2},y),(x+\frac{m}{2},y+\frac{n}{2})\right\}.$$
%We assume, without loss of generality, that $G=\{D(x,y):(x,y)\in \mathbb{Z}_\frac{m}{2}\times \mathbb{Z}_\frac{n}{2}\}$. The group operation in $G$ is defined by $D(x,y)\pm D(x',y')=D((x,y)\pm(x',y'))$, where the latter ``$\pm$" performs in $\mathbb{Z}_\frac{m}{2}\times \mathbb{Z}_\frac{n}{2}$.
The following proposition is straightforward from group theory.

\begin{Proposition}\label{1-pre-con-9}
\begin{enumerate}
\item[{\rm (1)}] For any distinct $(x,y)$ and $(x',y')$ from $\mathbb{Z}_\frac{m}{2}\times \mathbb{Z}_\frac{n}{2}$, $D(x,y)\cap D(x',y')=\emptyset.$
\item[{\rm (2)}] $\bigcup \limits_{(x,y)\in \mathbb{Z}_\frac{m}{2}\times \mathbb{Z}_\frac{n}{2}} D(x,y)=\mathbb{Z}_m\times \mathbb{Z}_n.$
\end{enumerate}
\end{Proposition}

Recall that $A=\mathbb Z_m\times \mathbb Z_n$ and
\begin{center}
$\begin{array}{ll}
A_{ee}=\left\{(x,y)\in A:\ x,y\equiv 0\ ({\rm mod }\ 2)\right\}, &
A_{eo}=\left\{(x,y)\in A:\ x\equiv 0\ ({\rm mod }\ 2), y\equiv 1\ ({\rm mod }\ 2)\right\}, \\[0.5em]
A_{oo}=\left\{(x,y)\in A:\ x,y\equiv 1\ ({\rm mod }\ 2)\right\}, &
A_{oe}=\left\{(x,y)\in A:\ x\equiv 1\ ({\rm mod }\ 2), y\equiv 0\ ({\rm mod }\ 2)\right\}.\\[0.5em]
\end{array}$
\end{center}

\begin{Proposition}\label{2-pre-con-9}
Let $m,n\equiv 2\ ({\rm mod}\ 4)$. For any $(x,y)\in \mathbb{Z}_\frac{m}{2}\times \mathbb{Z}_\frac{n}{2}$,
$$|D(x,y)\cap A_{oo}|=|D(x,y)\cap A_{oe}|=|D(x,y)\cap A_{eo}|=|D(x,y)\cap A_{ee}|=1.$$
\end{Proposition}
%\proof Note that $\frac{m}{2},\frac{n}{2}\equiv 1\ ({\rm mod}\ 2)$, the conclusion follows by checking the parity of the elements of $D(x,y)$ in $\mathbb{Z}_m\times \mathbb{Z}_n$. \qed

\begin{Construction}\label{2-mul}{\rm{(Doubling Construction)}} Let $m,n\equiv 2\ ({\rm mod}\ 4)$.
Suppose there exists an $(\frac{m}{2},\frac{n}{2},3,1)$-{\rm OOSPC} $\cal F$ whose leave is $L$. Then there exists an $(m,n,3,2,1)$-{\rm OOSPC} with $5\left|{\cal F}\right|$ codewords whose leave is
\begin{eqnarray}\label{eqn-L}
L'=\left(\bigcup\limits_{(x,y)\in L\cup\{(0,0)\}}(D(x,y)\setminus A_{ee})\right)\cup\left(\bigcup\limits_{(x,y)\in L}\{(2x,2y)\}\right).
\end{eqnarray}
Especially, if the given $(\frac{m}{2},\frac{n}{2},3,1)$-{\rm OOSPC} is $(s, t)$-regular, then the resulting $(m,n,3,2,1)$-{\rm OOSPC} is $(2s, 2t)$-regular.
\end{Construction}

\proof For each codeword $F=\{(0,0),(x_1,y_1),(x_2,y_2)\}\in{\cal F}$, construct a set ${\cal B}_F$ which consists of the following five $3$-subsets in $A$:
\begin{center}
\begin{tabular}{llll}
 & $\{(0,0),\alpha_1,2\alpha_1\},$ & $\{(0,0),\alpha_2,2\alpha_2\},$ & $\{(0,0),\alpha_3,2\alpha_3\},$\\
 & $\{(0,0),\beta_1,\beta_2\},$ & $\{(0,0),\beta_3,\beta_4\},$ &\\
\end{tabular}
\end{center}
satisfying that $\{\alpha_1,\beta_1,\beta_3\}=D(x_1,y_1)\setminus A_{ee}$, $\{\alpha_2,\beta_2,\beta_4\}=D(x_2,y_2)\setminus A_{ee}$ and $\{\alpha_3,\beta_2-\beta_1,\beta_4-\beta_3\}=D(x_2-x_1,y_2-y_1)\setminus A_{ee}$. This can be done because for any $(x,y)\in \mathbb{Z}_\frac{m}{2}\times \mathbb{Z}_\frac{n}{2}$, by Proposition \ref{2-pre-con-9}, $|D(x,y)\cap A_{oo}|=|D(x,y)\cap A_{oe}|=|D(x,y)\cap A_{eo}|=|D(x,y)\cap A_{ee}|=1$. So we can take $\alpha_1\in A_{oo}$, $\beta_1\in A_{eo}$, $\beta_3\in A_{oe}$, $\alpha_2\in A_{eo}$, $\beta_2\in A_{oe}$, $\beta_4\in A_{oo}$ and $\alpha_3\in A_{oe}$.

Note that
$\pm2\alpha_1=\pm(2x_1,2y_1),$ $\pm2\alpha_2=\pm(2x_2,2y_2)$ and $\pm2\alpha_3=\pm(2(x_2-x_1),2(y_2-y_1))$. They are distinct elements in $A_{ee}$. It follows that
\begin{eqnarray*}
\Delta{\cal B}_F&=&(\pm\{\alpha_1,\beta_1,\beta_3\})\cup(\pm\{\alpha_2,\beta_2,\beta_4\})\cup
(\pm\{\alpha_3,\beta_2-\beta_1,\beta_4-\beta_3\})\cup(\pm\{2\alpha_1,2\alpha_2,2\alpha_3\})\\
&=&\left(\bigcup \limits_{(x,y)\in \Delta F} \left(D(x,y)\setminus A_{ee}\right)\right)\cup\left(\bigcup \limits_{(x,y)\in \Delta F} \{(2x,2y)\}\right).
\end{eqnarray*}
Let ${\cal B}=\bigcup \limits_{F\in{\cal F}} {\cal B}_F$ and $V=\mathbb Z_{\frac{m}{2}}\times \mathbb Z_{\frac{n}{2}}\setminus (L\cup\{(0,0)\})$. Then
\begin{eqnarray*}
\Delta{\cal B}=\left(\bigcup\limits_{(x,y)\in V}(D(x,y)\setminus A_{ee})\right)\cup\left(\bigcup\limits_{(x,y)\in V}\{(2x,2y)\}\right).
\end{eqnarray*}
Thus ${\cal B}$ forms an $(m,n,3,2,1)$-OOSPC with $5\left|{\cal F}\right|$ codewords whose leave is of the form \eqref{eqn-L}.

Especially, if the given $(\frac{m}{2},\frac{n}{2},3,1)$-OOSPC is $(s, t)$-regular, then its leave $L$ along with $\{(0,0)\}$ forms an additive subgroup $S\times T$ of $\mathbb{Z}_{\frac{m}{2}}\times \mathbb{Z}_{\frac{n}{2}}$, where $S$ and $T$ are, respectively, the additive subgroups of order $s$ in $\mathbb{Z}_{\frac{m}{2}}$ and order $t$ in $\mathbb{Z}_{\frac{n}{2}}$. It is readily checked that the leave $L'$ of the resulting $(m,n,3,2,1)$-OOSPC satisfies
\begin{eqnarray*}
L'\cup\{(0,0)\}&=&\left(\bigcup\limits_{(x,y)\in S\times T}(D(x,y)\setminus A_{ee})\right)\cup\left(\bigcup\limits_{(x,y)\in S\times T}\{(2x,2y)\}\right)\\
&=&\left(\bigcup\limits_{(x,y)\in S\times T}D(x,y)\right)=S'\times T',
\end{eqnarray*}
where $S'$ and $T'$ are, respectively, the additive subgroups of order $2s$ in $\mathbb{Z}_{m}$ and order $2t$ in $\mathbb{Z}_{n}$.
Therefore, the resulting {\rm OOSPC} is $(2s, 2t)$-regular. \qed

\section{\label{4} Determination of $\Theta(m,n,3,\lambda_a,1)$ with $m,n\equiv 2\ ({\rm mod }\ 4)$}

This section is devoted to constructing optimal $(m,n,3,\lambda_a,1)$-OOSPCs with $\lambda_a=2,3$ for $m,n\equiv 2\ ({\rm mod }\ 4)$. In this case, $mn\equiv 4\ ({\rm mod }\ 8)$ and $\gcd(m,n,2)=2$. By Theorem \ref{c-main-1}, we have the following corollary.

\begin{Corollary}\label{cor-1}
For any $m,n\equiv 2\ ({\rm mod}\ 4)$ and $\lambda_a\in\{2,3\}$,
$$
\Theta(m,n,3,\lambda_a,1)\leq \left\lfloor\frac{5mn+4+8\omega}{24}\right\rfloor.
$$
\end{Corollary}

\begin{Proposition}\label{prop-equiv}
For any $m$ and $n$ such that $\gcd(mn,3)=1$, an $(m,n,3,2,1)$-{\rm OOSPC} is equivalent to an $(m,n,3,3,1)$-{\rm OOSPC}.
\end{Proposition}

\proof Let $\gcd(mn,3)=1$ and $X$ be a $3$-subset of $A$. By Lemma \ref{0-1}, $|{\rm{supp}}(\Delta X)|\geq 3$. Then by Lemma \ref{0-2}, $\lambda(X)\leq 2$. Hence, by the auto-correlation property $(1'')$, for $\gcd(mn,3)=1$, an $(m,n,3,2,1)$-{\rm OOSPC} is equivalent to an $(m,n,3,3,1)$-{\rm OOSPC}. \qed

\subsection{$(s, t)$-regular $(m,n,3,\lambda_a,1)$-OOSPCs}

\begin{Lemma}{\rm \cite[Theorem 4.8]{PC3}}\label{mn1}
For any $m$ and $n$ such that $mn\equiv1\ ({\rm mod}\ 6)$, there exists a $(1, 1)$-regular $(m,n,3,1)$-{\rm OOSPC}.
\end{Lemma}

\begin{Lemma}\label{regular-m5n3}
There exists a $(1, 3)$-regular $(m,n,3,1)$-{\rm OOSPC} for any $m\equiv 1,5\ ({\rm mod}\ 6)$ and $n\equiv3\ ({\rm mod}\ 6)$ except for $(m,n)=(1,9)$.
\end{Lemma}

\proof Let $m\equiv 1,5\ ({\rm mod}\ 6)$. For $n\in\{3,9\}$, due to $\gcd(m,3)=\gcd(m,9)=1$, a $(1, 3)$-regular $(m,n,3,1)$-OOSPC is equivalent to a cyclic Steiner triple system (CSTS) of order $mn$. It is known that a CSTS$(mn)$ exists if and only if $mn\equiv1,3\ ({\rm mod}\ 6)$ and $mn\neq 9$ (see \cite[Theorem 2.25]{CJC2}).

For $n\geq15$, start from a $(1, 3)$-regular $(1,n,3,1)$-OOSPC, which is equivalent to a CSTS$(n)$. Apply Construction \ref{R-DM} with an $(m,3,1)$-CDM to obtain an $(m, 3)$-regular $(m,n,3,1)$-{\rm OOSPC}. Then apply Construction \ref{Hole} to obtain a $(1, 3)$-regular $(m,n,3,1)$-{\rm OOSPC}. \qed

\begin{Lemma}\label{regular-add}
There exists a $(1, 3)$-regular $(m,n,3,2,1)$-{\rm OOSPC} with $\frac{5mn-12}{24}$ codewords for any $m\equiv 2,10\ ({\rm mod}\ 12)$ and $n\equiv 6\ ({\rm mod}\ 12)$.
\end{Lemma}

\proof For $(m,n)\in\{(2,6),(2, 18)\}$, all $\frac{5mn-12}{24}$ codewords of a $(1, 3)$-regular $(m,n,3,2,1)$-OOSPC are listed below:
\begin{center}
\begin{tabular}{llll}
$(m,n)=(2,6)$: & $\{(0,0),(0,1),(1,2)\}$, &  $\{(0,0),(0,3),(1,0)\}$; & \\
$(m,n)=(2,18)$: & $\{(0,0),(0,1),(0,2)\}$, &  $\{(0,0),(1,4),(0,8)\}$, & $\{(0,0),(1,7),(0,14)\}$,\\
			 & $\{(0,0),(0,3),(1,1)\}$, & $\{(0,0),(0,5),(1,8)\}$, & $\{(0,0),(0,7),(1,12)\}$, \\
			 & $\{(0,0),(0,9),(1,0)\}$.
\end{tabular}
\end{center}

For $m\equiv 2,10\ ({\rm mod}\ 12)$, $n\equiv 6\ ({\rm mod}\ 12)$ and $(m,n)\neq(2,18)$, by Lemma \ref{regular-m5n3}, there exists a $(1, 3)$-regular $(\frac{m}{2},\frac{n}{2},3,1)$-OOSPC with $\frac{mn-12}{24}$ codewords. Apply Construction \ref{2-mul} to obtain a $(2, 6)$-regular $(m,n,3,2,1)$-OOSPC with $\frac{5mn-60}{24}$ codewords. Then apply Construction \ref{Hole} with a $(1, 3)$-regular $(2,6,3,2,1)$-OOSPC, which has $2$ codewords, to obtain a $(1, 3)$-regular $(m,n,3,2,1)$-OOSPC with $\frac{5mn-12}{24}$ codewords. \qed

\begin{Lemma}\label{regular-mn3}
There exists a $(3, 3)$-regular $(m,n,3,1)${\rm-OOSPC} for any $m,n\equiv3\ ({\rm mod}\ 6)$, $m\neq 9$ and $n\neq 9$.
\end{Lemma}

\proof Without loss of generality, assume that $n\geq m$.

\textbf{Case 1:} $m=3$. When $n=3$, a $(3, 3)$-regular $(3,3,3,1)${\rm-OOSPC} is trivial.
When $n>9$, start from a $(1, 3)$-regular $(1,n,3,1)$-OOSPC, which is equivalent to a CSTS$(n)$.
Then apply Construction \ref{R-DM} with a $(3,3,1)$-CDM to obtain a $(3, 3)$-regular $(3,n,3,1)$-OOSPC.

\textbf{Case 2:} $m>9$. Then $n>9$. Start from a $(3, 3)$-regular $(3,n,3,1)$-OOSPC, which exists by Case 1, and apply Construction \ref{R-DM} with an $(\frac{m}{3},3,1)$-CDM to obtain an $(m, 3)$-regular $(m,n,3,1)$-OOSPC. Then apply Construction \ref{Hole} with a $(3, 3)$-regular $(m,3,3,1)$-OOSPC to obtain a $(3, 3)$-regular $(m,n,3,1)$-OOSPC. \qed

\begin{Lemma}\label{regular-mn3-1}
There exists a $(3, 3)$-regular $(9,n,3,1)${\rm-OOSPC} for any $n\equiv3\ ({\rm mod}\ 6)$.
\end{Lemma}

\proof For $n=3$, a $(3, 3)$-regular $(9,3,3,1)$-OOSPC has three codewords:
$\{(0,0),(1,0),(2,1)\}$, $\{(0,0),(1,2),(5,0)\}$, $\{(0,0),(2,0),(4,2)\}$.
For $n=9$, there exists a $(3, 3)$-regular $(9,9,3,1)$-OOSPC by \cite[Lemma 4.5]{PC3}.
For $n\geq15$, start from the $(3, 3)$-regular $(9,3,3,1)$-OOSPC, and apply Construction \ref{R-DM} with an $(\frac{n}{3},3,1)$-CDM to obtain a $(3, n)$-regular $(9,n,3,1)$-OOSPC. Then apply Construction \ref{Hole} with a $(3, 3)$-regular $(3,n,3,1)$-OOSPC, which exists by Lemma \ref{regular-mn3}, to obtain a $(3, 3)$-regular $(9,n,3,1)$-OOSPC. \qed

Denote by $[a, b]$ the set of integers $v$ such that $a\leq v \leq b$.

\begin{Lemma} \label{regular-m2n2} For any $m\equiv 2\ ({\rm mod}\ 12)$ and $n\equiv 2\ ({\rm mod}\ 4)$,
there exists a $(2, n)$-regular $(m,n,3,2,1)$-{\rm OOSPC} with $\frac{5n(m-2)}{24}$ codewords.
\end{Lemma}

\proof For $m=2$, the conclusion is trivial. For $m\geq 14$, all $\frac{5n(m-2)}{24}$ codewords of a $(2, n)$-regular $(m,n,3,2,1)$-OOSPC are listed as follows:
\begin{center}
\begin{tabular}{ll}	
Type 4: & $\{(0,0),(12t+3,2i),(24t+6,4i)\},$  \\
            & $\{(0,0),(12t+5,2i+1),(24t+10,4i+2)\},$  \\
            & $\{(0,0),(12s+9,2i+1),(24s+18,4i+2)\},$  \\
            & $\{(0,0),(12s+11,2i),(24s+22,4i)\},$  \\
            & $\{(0,0),(\frac{m}{2}-6j-1,2i+1),(m-12j-2,4i+2)\},$ \\
Type 6: & $\{(0,0),(6j+1,2i),(12j+3,4i+1)\},$\\
            & $\{(0,0),(6j+4,2i+1),(12j+5,4i+2)\},$\\
\end{tabular}
\end{center}
where $i\in[0,\frac{n}{2}-1]$, $j\in[0,\frac{m-14}{12}]$, $s\in[0,\lfloor\frac{m-26}{24}\rfloor]$ and $t\in[0,\lfloor\frac{m-14}{24}\rfloor]$. \qed

\subsection{Proof of Theorem \ref{c-main-2}}

\begin{Lemma}\label{mn2mn10}
$\Theta(m,n,3,2,1)=\frac{5mn+4}{24}$ for any $m,n\equiv 2\ ({\rm mod}\ 12)$ and $m,n\equiv10\ ({\rm mod}\ 12)$.
\end{Lemma}

\proof $m,n\equiv 2\ ({\rm mod}\ 12)$ and $m,n\equiv10\ ({\rm mod}\ 12)$ both imply $mn\equiv 4\ ({\rm mod}\ 24)$. By Corollary \ref{cor-1}, $\Theta(m,n,3,2,1)\leq\frac{5mn+4}{24}$. For $(m,n)=(2,2)$, an optimal $(2,2,3,2,1)$-OOSPC has only one codeword $\{(0,0),(1,0),(0,1)\}$.
For $(m,n)\neq(2,2)$, start from a $(1, 1)$-regular $(\frac{m}{2},\frac{n}{2},3,1)$-OOSPC with $\frac{mn-4}{24}$ codewords, which exists by Lemma \ref{mn1}, and apply Construction \ref{2-mul} to obtain a $(2, 2)$-regular $(m,n,3,2,1)$-{\rm OOSPC} with $\frac{5mn-20}{24}$ codewords. Then apply Construction \ref{Hole} with an optimal $(2,2,3,2,1)$-OOSPC with $1$ codeword to obtain an optimal $(m,n,3,2,1)$-OOSPC with $\frac{5mn+4}{24}$ codewords. \qed

\begin{Lemma}\label{m2n10}
$\Theta(m,n,3,2,1)=\frac{5mn-4}{24}$ for any $m\equiv 2\ ({\rm mod}\ 12)$ and $n\equiv 10\ ({\rm mod}\ 12)$.
\end{Lemma}

\proof By Corollary \ref{cor-1}, $\Theta(m,n,3,2,1)\leq\frac{5mn-4}{24}$.
For $m=2$, set $n=2t$, where $t\equiv 5\ ({\rm mod}\ 6)$.
We here give an explicit construction for an optimal $(2,n,3,2,1)$-{OOSPC} on $\mathbb{Z}_2\times \mathbb{Z}_2\times \mathbb{Z}_{t}$ with $\frac{5n-2}{12}=\frac{5t-1}{6}$ codewords:
\begin{center}
\noindent\begin{tabular}{lll}
Type 3: & $\{(0,0,0),(0,1,0),(1,0,0)\},$ & \\[0.2em]
Type 4: & $\{(0,0,0),(0,1,2i-1),(0,0,4i-2)\}$,&  $i\in[1, \frac{t+1}{6}]$,\\[0.2em]
             & $\{(0,0,0),(1,0,2i),(0,0,4i)\}$,&  $i\in[1, \frac{t+1}{6}]$,\\[0.2em]
             & $\{(0,0,0),(1,1,\frac{t+1}{2}-i),(0,0,t-2i+1)\}$,& $i\in[1, \frac{t-5}{6}]$,\\[0.2em]
Type 6: & $\{(0,0,0),(1,1,2i-1),(0,1,4i-2)\}$,& $i\in[1, \frac{t+1}{6}]$,\\[0.2em]
             & $\{(0,0,0),(0,1,4i),(1,0,\frac{t+1}{3}+2i)\}$,& $i\in[1, \frac{t-5}{6}]$.
\end{tabular}
\end{center}
%Note that ${\rm DL}({\cal F})=\{\pm(1,0, \frac{t+1}{3})\}$.

For $m\geq 14$, begin with a $(2, n)$-regular $(m,n,3,2,1)$-{\rm OOSPC} with $\frac{5n(m-2)}{24}$ codewords, which comes from Lemma \ref{regular-m2n2}. Apply Construction \ref{Hole} with a $(2,n,3,2,1)$-OOSPC with $\frac{5n-2}{12}$ codewords to obtain an optimal $(m,n,3,2,1)$-{\rm OOSPC} with $\frac{5mn-4}{24}$ codewords. \qed

\begin{Lemma}\label{m10n6}
For any $m\equiv 2,10\ ({\rm mod}\ 12)$ and $n\equiv 6\ ({\rm mod}\ 12)$, $\Theta(m,n,3,2,1)=\frac{5mn-12}{24}$ and $\Theta(m,n,3,3,1)=\frac{5mn+12}{24}$.
\end{Lemma}

\proof The condition $m\equiv 2,10\ ({\rm mod}\ 12)$ and $n\equiv 6\ ({\rm mod}\ 12)$ implies $mn\equiv 12\ ({\rm mod}\ 24)$. By Corollary \ref{cor-1}, $\Theta(m,n,3,2,1)\leq\frac{5mn-12}{24}$ and $\Theta(m,n,3,3,1)\leq\frac{5mn+12}{24}$. By Lemma \ref{regular-add},
there exists a $(1, 3)$-regular $(m,n,3,2,1)$-{\rm OOSPC} with $\frac{5mn-12}{24}$ codewords for any $m\equiv 2,10\ ({\rm mod}\ 12)$ and $n\equiv 6\ ({\rm mod}\ 12)$. Thus $\Theta(m,n,3,2,1)=\frac{5mn-12}{24}$.

Start from the resulting $(1, 3)$-regular OOSPC, and apply Construction \ref{Hole} with an optimal $(1,3,3,3,1)$-OOSPC, which consists of the unique codeword $\{(0,0),(0,1),(0,2)\}$, to obtain an optimal $(m,n,3,3,1)$-OOSPC with $\frac{5mn+12}{24}$ codewords. Thus $\Theta(m,n,3,3,1)=\frac{5mn+12}{24}$. \qed

\begin{Lemma}\label{mn6}
For any $m,n\equiv 6\ ({\rm mod}\ 12)$, $\Theta(m,n,3,2,1)=\frac{5mn-12}{24}$ and $\Theta(m,n,3,3,1)=\frac{5mn+36}{24}$.
\end{Lemma}

\proof The condition $m,n\equiv 6\ ({\rm mod}\ 12)$ implies $mn\equiv 12\ ({\rm mod}\ 24)$. By Corollary \ref{cor-1}, $\Theta(m,n,3,2,1)\leq\frac{5mn-12}{24}$ and $\Theta(m,n,3,3,1)\leq\frac{5mn+36}{24}$. When $(m,n)=(6,6)$, the conclusion follows from Example \ref{eg}. When $(m,n)\neq(6,6)$, there is a $(3, 3)$-regular $(\frac{m}{2},\frac{n}{2},3,1)$-OOSPC with $\frac{mn-36}{24}$ codewords by Lemmas \ref{regular-mn3} and \ref{regular-mn3-1}. Apply Construction \ref{2-mul} to obtain a $(6, 6)$-regular $(m,n,3,2,1)$-OOSPC with $\frac{5mn-180}{24}$ codewords. Then apply Construction \ref{Hole} with a $(6,6,3,2,1)$-OOSPC with 7 codewords to obtain a $(m,n,3,2,1)$-OOSPC with $\frac{5mn-12}{24}$ codewords. Thus $\Theta(m,n,3,2,1)=\frac{5mn-12}{24}$. Apply Construction \ref{Hole} with a $(6,6,3,3,1)$-OOSPC with 9 codewords to obtain a $(m,n,3,3,1)$-OOSPC with $\frac{5mn+36}{24}$ codewords. \qed

By Proposition \ref{prop-equiv}, for any $m$ and $n$ such that $\gcd(mn,3)=1$, an $(m,n,3,2,1)$-{\rm OOSPC} is equivalent to an $(m,n,3,3,1)$-{\rm OOSPC}. Note that $\Theta(m,n,3,\lambda_a,1)=\Theta(n,m,3,\lambda_a,1)$. Now combining the results from Lemmas \ref{mn2mn10}-\ref{mn6}, one can complete the proof of Theorem \ref{c-main-2}.

\section{\label{conclusion}Concluding remarks}

Compared with \eqref{eqq3}, Theorem \ref{c-main-1} provides a much more complicated upper bound on the size of an $(m,n,3,\lambda_a,1)$-OOSPC with $\lambda_a\in\{2,3\}$.
%When $mn\equiv 0\ ({\rm mod }\ 4)$, $\Theta(m,n,k,\lambda_a,\lambda_c)=\frac{mn-2}{4}$ by Theorem \ref{thm-sawa}.
It seems that this bound is good for $mn\equiv 0\ ({\rm mod }\ 4)$. On one hand, when $\gcd(m,n)=1$, an
$(m,n,k,\lambda_a,\lambda_c)$-OOSPC is equivalent to a $1$-D $(mn,k,\lambda_a,\lambda_c)$-OOC \cite{GW}. Let $\Phi(mn,k,\lambda_a,\lambda_c)$ denote the largest possible size among all $1$-D $(mn,k,\lambda_a,\lambda_c)$-OOCs. Then $\Theta(m,n,k,\lambda_a,\lambda_c)=\Phi(mn,k,\lambda_a,\lambda_c)$ for $\gcd(m,n)=1$. When $v\equiv 0\ ({\rm mod }\ 4)$, the exact value of $\Phi(v,3,\lambda_a,1)$ has been determined in the literature. Note that a $1$-D $(v,k,k,1)$-OOC is often referred to as a {\em conflict-avoiding code}, which finds its application on a multiple-access collision channel without feedback.

\begin{table}[t]
\begin{center}
\setlength{\belowcaptionskip}{10pt}	
\begin{tabular}{|c|c|c|c|}\hline	
\multirow{2}*{$m$} & \multirow{2}*{$n$}  & \multicolumn{2}{c|}{$\Theta(m,n,3,\lambda_a,1)$} \\\cline{3-4}
& & $\lambda_a=3$ & $\lambda_a=2$\\\hline
$2$ & $4$ & $1$ & $1$\\\hline
$2$ & $8$ & $3$ & $3$\\\hline
$2$ & $12$ & $ 5$ & $5$\\\hline
$2$ & $16$ & $7$ & $7$\\\hline
$2$ & $20$ & $8 $ & $8 $\\\hline
$2$ & $24$ & $10$ & $10$ \\\hline
$2$ & $28$ & $12 $ & $12 $\\\hline
$2$ & $32$ & $13$ & $13$ \\\hline
$2$ & $36$ & $15 $ & $15 $\\\hline
$2$ & $40$ & $16$ & $16$ \\\hline
$2$ & $44$ & $18 $ & $18 $\\\hline
$2$ & $48$ & $20$ & $20$ \\\hline
$2$ & $52$ & $21 $ & $21 $\\\hline
$2$ & $56$ & $23$ & $23$ \\\hline
$2$ & $60$ & $25 $ & $25 $\\\hline
$2$ & $64$ & $26$ & $26$ \\\hline
$2$ & $68$ & $28 $ & $28 $\\\hline
$2$ & $72$ & $30$ & $29$ \\\hline
\end{tabular}
\begin{tabular}{|c|c|c|c|}\hline	
\multirow{2}*{$m$} & \multirow{2}*{$n$}  & \multicolumn{2}{c|}{$\Theta(m,n,3,\lambda_a,1)$} \\\cline{3-4}
& & $\lambda_a=3$ & $\lambda_a=2$\\\hline
$3$ & $12$ & $9$ & $7$ \\\hline
$3$ & $24$ & $17$ & $15$ \\\hline
$3$ & $36$ & $25$ & $23$ \\\hline
$3$ & $48$ & $33$ & $31$ \\\hline
$4$ & $4$ & $4$ & $4$  \\\hline
$4$ & $6$ & $5$ & $ 5$ \\\hline
$4$ & $8$ & $6$ & $6$\\\hline
$4$ & $10$ & $ 8$ & $8 $\\\hline
$4$ & $12$ & $11$ & $10$ \\\hline
$4$ & $14$ & $12 $ & $12 $\\\hline
$4$ & $16$ & $13$ & $13$  \\\hline
$4$ & $18$ & $15 $ & $15 $\\\hline
$4$ & $20$ & $17$ & $17$\\\hline
$4$ & $22$ & $ 18$ & $18 $\\\hline
$4$ & $24$ & $20$ & $19$ \\\hline
$4$ & $26$ & $21 $ & $21 $\\\hline
$4$ & $28$ & $23$ & $23$\\\hline
$4$ & $30$ & $25 $ & $25 $\\\hline
\end{tabular}
\begin{tabular}{|c|c|c|c|}\hline	
\multirow{2}*{$m$} & \multirow{2}*{$n$}  & \multicolumn{2}{c|}{$\Theta(m,n,3,\lambda_a,1)$} \\\cline{3-4}
& & $\lambda_a=3$ & $\lambda_a=2$\\\hline
$4$ & $32$ & $26$ & $26$ \\\hline
$4$ & $34$ & $28 $ & $28 $\\\hline
$4$ & $36$ & $30$ & $29$  \\\hline
$5$ & $20$ & $22$ & $22$ \\\hline
$6$ & $8$ & $10$ & $10$ \\\hline
$6$ & $12$ & $16 $ & $15$\\\hline
$6$ & $16$ & $20$ & $20$ \\\hline
$6$ & $20$ & $25 $ & $25 $\\\hline
$6$ & $24$ & $31$ & $29$ \\\hline
$8$ & $8$ & $13$ & $13$ \\\hline
$8$ & $10$ & $16$ & $16$ \\\hline
$8$ & $12$ & $20$ & $19$\\\hline
$8$ & $14$ & $23$ & $23$ \\\hline
$8$ & $16$ & $25$ & $25$\\\hline
$8$ & $18$ & $30$ & $29$ \\\hline
$9$ & $12$ & $25$ & $23$  \\\hline
$10$ & $12$ & $25 $ & $25 $\\\hline
$12$ & $12$ & $32$ & $29$ \\\hline
\end{tabular}
\end{center}
\caption{$\Theta(m,n,3,\lambda_a,1)$ not covered by Theorem \ref{c-main-2} for $mn\equiv 0\ ({\rm mod }\ 4)$, $mn\leq 150$ and   $\gcd(m,n)\neq1$}
\end{table}

\begin{Theorem}\label{thm OOC-1}{\rm \cite{FWW,WaChF}}
\begin{eqnarray*}
\Phi(v,3,2,1)=\left\{
\begin{array}{ll}
\left\lfloor\frac{7v}{32}\right\rfloor,&\hbox{ \rm{if} $v\equiv 0\ ({\rm mod }\ 8)$ \rm{and} $v\neq 64$\rm{;}}\\[0.4em]
13, & \hbox{ \rm{if} $v=64$}\rm{;} \\[0.4em]
\left\lfloor\frac{7v+4}{32}\right\rfloor,&\hbox{ \rm{if} $v\equiv 4\ ({\rm mod }\ 8)$\rm{.}}\\
\end{array}
\right.
\end{eqnarray*}
\end{Theorem}

\begin{Theorem}{\rm \cite{LT,JMJ,MFU,FLM}}\label{cac}
\begin{eqnarray*}
\Phi(v,3,3,1)=\left\{
\begin{array}{ll}
\left\lfloor\frac{7v+16}{32}\right\rfloor,&\hbox{ \rm{if }$v\equiv 0\ ({\rm mod }\ 24)$ \rm{and} $v\neq 48$\rm{;}}\\[0.4em]
10, & \hbox{ \rm{if} $v=48$}; \\[0.4em]
\left\lfloor\frac{7v+4}{32}\right\rfloor,&\hbox{ \rm{if }$v\equiv 4,20\ ({\rm mod }\ 24)$\rm{;}}\\[0.4em]
\left\lfloor\frac{7v}{32}\right\rfloor,&\hbox{ \rm{if }$v\equiv 8,16\ ({\rm mod }\ 24)$ \rm{and} $v\neq 64$\rm{;}}\\[0.4em]
\ 13, & \hbox{ \rm{if} $v=64$};\\[0.4em]
\left\lfloor\frac{7v+20}{32}\right\rfloor,&\hbox{ \rm{if }$v\equiv 12\ ({\rm mod }\ 24)$\rm{.}}\\
\end{array}
\right.
\end{eqnarray*}
\end{Theorem}

It is easy to check that Theorems \ref{thm OOC-1} and \ref{cac} satisfy the bound for $\Theta(m,n,3,\lambda_a,1)$ in Theorem \ref{c-main-1} when $\gcd(m,n)=1$ except for $mn\in\{48,64\}$.

On the other hand, when $\gcd(m,n)\neq1$, Theorem \ref{c-main-2} determines the values of $\Theta(m,n,3,\lambda_a,1)$ with $\lambda_a=2,3$ for $m,n\equiv 2\ ({\rm mod }\ 4)$, which coincides with the bound in Theorem \ref{c-main-1}. By computer search, it is shown that for any $m$ and $n$ such that $mn\equiv 0\ ({\rm mod }\ 4)$ and $mn\leq 150$, there exists an $(m,n,3,\lambda_a,1)$-OOSPC attaining the bound in Theorem \ref{c-main-1} (see Table 1). The interested reader may get a copy of these data from the authors. We conjecture that when $mn\equiv 0\ ({\rm mod }\ 4)$, our bound for $\Theta(m,n,3,\lambda_a,1)$ with $\lambda_a\in\{2,3\}$ shown in Theorem \ref{c-main-1} is tight.

Theorem \ref{c-main-2} determines the value of $\Theta(m,n,3,\lambda_a,1)$ with $\lambda_a\in\{2,3\}$ for $m,n\equiv 2\ ({\rm mod}\ 4)$. To prove Theorem \ref{c-main-2}, the doubling construction (Construction \ref{2-mul}) plays an important role. It seems that to solve other cases of $m$ and $n$ such that $mn\equiv 0\ ({\rm mod }\ 4)$, one must explore a quadrupling construction.

%We mentioned again that in Section \ref{4} we succeeded to complete the $\Theta(m,n,3,2,1)$ for $m,n\equiv 2\ ({\rm mod }\ 4)$. We believe that with a similar technique, $\Theta(m,n,3,3,1)$ for $m,n\equiv 2\ ({\rm mod }\ 4)$ can also be solved. In this case, one only needs to consider the additional condition of $mn\equiv 0\ ({\rm mod }\ 3)$, since an $(m,n,3,2,1)$-OOSPC is actually an $(m,n,3,3,1)$-OOSPC if $mn\not\equiv 0\ ({\rm mod }\ 3)$. Besides, observing the combinatorial constructions utilized in Section \ref{4}, it seems that the structure of an optimal OOSCP is irregular, which makes it difficult to find effective recursive constructions for optimal $(m,n,3,\lambda_a,1)$-OOSCPs. By tedious computation and analysis, direct constructions certainly always work, but it could not be a good way. A better technique is expected for the problem.


\begin{thebibliography}{Z}
\bibitem{JAM}
R.J.R. Abel and M. Buratti, Some progress on $(v,4,1)$ difference families and optical orthogonal codes,
{\it J. Combin. Theory Ser. A}, Vol. 106 (2004), 59--75.

\bibitem{BM}
M. Buratti, Cyclic designs with block size $4$ and related optimal optical orthogonal codes,
{\it Des. Codes Cryptography}, Vol. 26 (2002), 111--125.

\bibitem{BP}
M. Buratti and A. Pasotti, Further progress on difference families with block size $4$ or $5$,
{\it Des. Codes Cryptography}, Vol. 56 (2010), 1--20.

\bibitem{BPW}
M. Buratti, A. Pasotti, D. Wu, On optimal $(v,5,2,1)$ optical orthogonal codes,
{\it Des. Codes Cryptography}, Vol. 68 (2013), 349--371.

%\bibitem{CFM}
%Y. Chang, R. Fuji-Hara, and Y. Miao, {\it Combinatorial constructions of optimal optical orthogonal codes with weight $4$,} IEEE Trans. Inform. Theory, Vol. 49, No. 5 (2003), 1283--1292.

%\bibitem{CJ}
%Y. Chang and L. Ji, {\it Optimal $(4up,5,1)$ optical orthogonal codes}, J. Combin. Des., Vol. 12, No. 5 (2004), 346--361.

%\bibitem{ChM}
%Y. Chang and Y. Miao, {\it Constructions for optimal optical orthogonal codes}, Discrete Math., Vol. 261 (2003), 127--139.

%\bibitem{ChY}
%Y. Chang and J. Yin, Further results on optimal optical orthogonal codes with weight $4$, {\it Discrete Math.}, Vol. 279 (2004), 135--151.

\bibitem{CJL}
J. Chen, L. Ji, and Y. Li, New optical orthogonal signature pattern codes with maximum collision parameter $2$ and weight $4$, {\it Des. Codes Cryptography}, Vol. 85 (2017), 299--318.

\bibitem{CJL2}
J. Chen, L. Ji, and Y. Li, Combinatorial constructions of optimal $(m, n, 4, 2)$ optical orthogonal signature pattern codes, {\it Des. Codes Cryptography}, Vol. 86 (2018), 1499--1525.

\bibitem{CSW}
F.R.K. Chung, J.A Salehi, and V.K. Wei, Optical orthogonal codes: design, analysis, and applications, {\it IEEE Trans. Inform. Theory}, Vol. 35, No. 3 (1989), 595--604.

\bibitem{CJC}
C.J. Colbourn, {\it Difference matrices},
in: CRC Handbook of Combinatorial Designs (C.J. Colbourn and J.H. Dinitz, eds.),
CRC Press, Boca Raton, 2007, 411--419.

\bibitem{CJC2}
C.J. Colbourn, {\it Triple systems},
in: CRC Handbook of Combinatorial Designs (C.J. Colbourn and J.H. Dinitz, eds.),
CRC Press, Boca Raton, 2007, 58--62.

%\bibitem{FCJ}
%T. Feng, Y. Chang, and L. Ji, {\it Constructions for rotatinal Steiner quadruple systems,}
%J. Combin. Des., Vol. 17, No. 5 (2009), 353--368.

\bibitem{FCJ2}
T. Feng, Y. Chang, and L. Ji, Constructions for strictly cyclic $3$-designs and applications to optimal OOCs with $\lambda=2$, {\it J. Combin. Theory Ser. A}, Vol. 115 (2008), 1527--1551.

\bibitem{FWW}
T. Feng, L. Wang, and X. Wang, Optimal $2$-D $(n \times m, 3, 2, 1)$-optical orthogonal codes and related equi-difference conflict avoiding codes, {\it Des. Codes Cryptography}, Vol. 87 (2019), 1499--1520.

\bibitem{FLM}
H. Fu, Y. Lin, and M. Mishima, Optimal conflict-avoiding codes of even length and weight $3$,
{\it IEEE Trans. Inform. Theory}, Vol. 56, No. 11 (2010), 5747--5756.

\bibitem{FM}
R. Fuji-Hara and Y. Miao, Optical orthogonal codes: their bounds and new optimal constructions,
{\it IEEE Trans. Inform. Theory}, Vol. 46, No. 7 (2000), 2396--2406.

%\bibitem{FMY}
%R. Fuji-Hara, Y. Miao, and J. Yin, {\it Optimal $(9v,4,1)$ optical orthogonal codes,}
%SIAM J. Discrete Math., Vol. 14 (2001), 256--266.

%\bibitem{GMY}
%G. Ge, Y. Miao, and X. Sun, {\it Perfect difference families, perfect difference matrices, and related combinatorial structures},
%J. Combin. Designs, Vol. 18, No. 6 (2010), 415--449.

\bibitem{GY}
G. Ge and J. Yin, Constructions for optimal $(v,4,1)$ optical orthogonal codes,
{\it IEEE Trans. Inform. Theory}, Vol. 47, No. 7 (2001), 2998--3004.

\bibitem{JDWG}
L. Ji, B. Ding, X. Wang, and G. Ge, Asymptotically optimal optical orthogonal signature pattern codes, {\it IEEE Trans. Inform. Theory}, Vol. 64, No. 7 (2018), 5419--5431.

\bibitem{JMJ}
M. Jimbo, M. Mishima, S. Janiszewski, A.Y. Teymorian, and V. D. Tonchev, On conflict-avoiding codes of length  $n=4m$ for three active users, {\it IEEE Trans. Inform. Theory}, Vol. 53, No. 8 (2007), 2732--2742.

\bibitem{J}
S.M. Johnson, A new upper bound for error-correcting codes, {\it IEEE Trans. Inform. Theory}, Vol. 8 (1962), 203--207.

\bibitem{KK}
K. Kitayama, Novel spatial spread spectrum based fiber optic CDMA networks for image transmission,
{\it IEEE J. Sel. Areas Commun.}, Vol. 12, No. 4 (1994), 762--772.

\bibitem{k}
K. Kitayama, Optical Code Division Multiple Access: A Practical Perspective, Cambridge University
Press, New York, 2014.

\bibitem{KY}
W.C. Kwong and G.C. Yang, Image transmission in multicore-fiber code-division multiple-access networks,
{\it IEEE Commun. Letters}, Vol. 2, No. 10 (1998), 285--287.

%\bibitem{KY2}
%W. C. Kwong and G. C. Yang, {\it Double-weight signature pattern codes for multicore-fiber code-division multiple-access networks}, IEEE Commun. Letters, Vol. 5, No. 5 (2001), 203--205.

\bibitem{LT}
V.I. Levenshtein and V.D. Tonchev, Optimal conflict-avoiding codes for three active users,
in \textit{Proc. IEEE Int. Symp. Inf. Theory}, (2005), 535--537.

%\bibitem{MC}
%S. Ma and Y. Chang, {\it A new class of optimal optical orthogonal codes with weight five,}
%IEEE Trans. Inform. Theory, Vol. 50, No. 8 (2004), 1848--1850.

%\bibitem{MC2}
%S. Ma and Y. Chang, {\it Constructions of optimal optical orthogonal codes with weight five,}
%J. Combin. Des., Vol. 13, No. 1 (2005), 54--69.

\bibitem{MFU}
M. Mishima, H. Fu, and S. Uruno, Optimal conflict-avoiding codes of length $n\equiv 0\ ({\rm mod }\ 16)$ and weight $3$, {\it Des. Codes Cryptography}, Vol. 52 (2009), 275--291.

\bibitem{PC}
R. Pan and Y. Chang, Combinatorial constructions for maximum optical orthogonal signature pattern codes,
{\it Discrete Math.}, Vol. 313 (2013), 2918--2931.

%\bibitem{PC2}
%R. Pan and Y. Chang, {\it Further results on optimal $(m,n,4,1)$ optical orthogonal signature pattern codes \emph{(in Chinese)},} Sci. Sin. Math., Vol. 44 (2014), 1141--1152.

\bibitem{PC3}
R. Pan and Y. Chang, $(m,n,3,1)$ optical orthogonal signature pattern codes with maximum possible size,
{\it IEEE Trans. Inform. Theory}, Vol. 61, No. 2 (2015), 1139--1148.

\bibitem{MS}
M. Sawa, Optical orthogonal signature pattern codes with maximum collision parameter $2$ and weight $4$,
{\it IEEE Trans. Inform. Theory}, Vol. 56, No. 7 (2010), 3613--3620.

\bibitem{sawa}
M. Sawa and S. Kageyama, Optimal optical orthogonal signature pattern codes of weight $3$,
{\it Biometrical Letters}, Vol. 46 (2009), 89--102.

%\bibitem{WaCh}
%X. Wang and Y. Chang, {\it Further results on $(v,4,1)$-perfect difference families,}
%Discrete Math., Vol. 310 (2010), 1995--2006.

%\bibitem{WaCh2}
%X. Wang and Y. Chang, {\it Further results on optimal $(v,4,2,1)$-OOCs,}
%Discrete Math., Vol. 312, No. 2 (2012), 331--340.

\bibitem{WaChF}
X. Wang, Y. Chang, and T. Feng, Optimal $2$-D $(n\times m,3,2,1)$-optical orthogonal codes,
{\it IEEE Trans. Inform. Theory}, Vol. 59, No. 1 (2013), 710--725.

\bibitem{GW}
G.C. Yang and W.C. Kwong, Two-dimensional spatial signature patterns,
{\it IEEE Trans. Commun.}, Vol. 44, No. 2 (1996), 184--191.

\bibitem{Y}
J. Yin, Some combinatorial constructions for optical orthogonal codes,
{\it Discrete Math.}, Vol. 185 (1998), 201--219.

%\bibitem{Y1}
%J. Yin, A general construction for optimal cyclic packing designs,
%{\it J. Combin. Theory Ser. A}, Vol. 97 (2002), 272--284.

%\bibitem{ZQ}
%H. Zhao and R. Qin, {\it Combinatorial constructions for optimal optical orthogonal signature pattern codes,}
%Discrete Math., Vol. 339 (2016), 179--193.
\end{thebibliography}
\end{document}